\documentclass[12pt,a4paper,reqno]{amsart}

\usepackage{amsfonts}
\usepackage{amsmath}
\usepackage{amssymb}
\usepackage{amsthm}
\usepackage{amsxtra}
\usepackage{bbm}
\usepackage{braket}
\usepackage{comment}
\usepackage{extarrows}
\usepackage{graphicx}
\usepackage{latexsym}
\usepackage{mathrsfs}
\usepackage{yhmath}
\usepackage{nicematrix}
\usepackage{mathtools}

\usepackage{float}
\usepackage{booktabs}
\usepackage{verbatim}
\usepackage{xcolor}

\usepackage[pdfborder={0 0 0}]{hyperref} 
\hypersetup{breaklinks, raiselinks}

\usepackage{bm}

\usepackage{cleveref}

\usepackage{a4wide}
\usepackage{fullpage}

\theoremstyle{plain}
\newtheorem{Theorem}{Theorem}[section]
\newtheorem{Lemma}[Theorem]{Lemma}
\newtheorem{Corollary}[Theorem]{Corollary}
\newtheorem{Proposition}[Theorem]{Proposition}

\theoremstyle{definition}

\newtheorem{Assumptions and Discussion}[Theorem]{Assumptions and Discussion}

\theoremstyle{remark}

\newtheorem*{acknowledgment*}{Acknowledgment}

\numberwithin{equation}{section}

\usepackage[shortlabels]{enumitem}
\SetEnumerateShortLabel{a}{\textup{(\alph*)}}
\SetEnumerateShortLabel{A}{\textup{(\Alph*)}}
\SetEnumerateShortLabel{1}{\textup{(\arabic*)}}
\SetEnumerateShortLabel{i}{\textup{(\roman*)}}
\SetEnumerateShortLabel{I}{\textup{(\Roman*)}}

\def\bar#1{\overline{#1}}

\begin{document}

\title{Compactness of conformal Chern-minimal surfaces in Hermitian surface}

\author[ Xiaowei Xu]{Xiaowei Xu}





\address{Xiaowei Xu, School of Mathematical Sciences, University of Science and Technology of China, Hefei, Anhui, 230026, P.R.~China; CAS, Wu Wen-Tsun Key Laboratory of Mathematics.}
\email{xwxu09@ustc.edu.cn}

\begin{abstract}
The Chern-minimal surfaces in Hermitian surface play a similar role as minimal surfaces in K\"ahler surface (see \cite{[PX-21]}) from the viewpoint of submanifolds.
This paper studies the compactness of Chern-minimal surfaces. We prove that any sequence $\{f_n\}$ of conformal Chern-minimal maps from closed Riemann surface
$(\Sigma,\emph{\texttt{j}})$ into a compact Hermitian surface $(M, J, h)$ with bounded area has a bubble tree limit, which consisting of a Chern-minimal map 
$f_0$ from $\Sigma$ into $M$ and a finite set of Chern-minimal maps from $S^2$ into $M$. We also show that the limit preserves area and homotopy class.

\end{abstract}
\maketitle

\section{Introduction}

It is known that there are abundant results on minimal surfaces in four-dimensional Riemannian manifold, especially the ones in K\"ahler surface. 
For instance, S. Webster (see \cite{[W-84]}, \cite{[W-86]}) counted the complex and anti-complex points of generic minimal immersion from Riemann surface into K\"ahler surface by using the topological data. 
Here the generic means that the immersion is neither holomorphic nor anti-holomorphic. More explicitly, let $f$ be a conformal minimal immersion from closed
Riemann surface $(\Sigma, \emph{\texttt{j}})$ into K\"ahler surface $(M, J, \omega)$, a point $x\in\Sigma$ is called \emph{complex} 
(resp. \emph{anti-complex}) if $df\circ\emph{\texttt{j}}=J\circ df$ (resp. $df\circ\emph{\texttt{j}}=-J\circ df$) holds at $x$, where $\emph{\texttt{j}}$ and $J$ are the complex structures.
S. Webster proved that the sum of complex points and 
anti-complex points of a generic minimal immersion $f$, denoted by $P$ and $Q$, respectively, can be counted by
\begin{equation}\label{Z-1}
P+Q=-(\chi(T\Sigma)+\chi(T^\perp\Sigma)),
\end{equation}
\begin{equation}\label{Z-2}
P-Q=-f^*c_1(M)[\Sigma],
\end{equation}
where $\chi(T\Sigma)$, $\chi(T^\perp\Sigma)$ are the Euler characteristic of tangent bundle $T\Sigma$, normal bundle $T^\perp\Sigma$, respectively, 
$c_1(M)$ is the first Chern class of $M$ and $[\Sigma]$ is the fundamental class of $\Sigma$. There are some progresses on Webster's
formulae. J.G. Wolfson (see \cite{[W-89]}) gave a new proof and deep applications of these formulae in the theory of minimal surfaces in K$\ddot{a}$hler surface;
J.Y. Chen and  G. Tian (see \cite{[CT-97]}), X.L. Han and J.Y. Li (see \cite{[HL-10]}) proved that (\ref{Z-1}) holds for minimal surfaces in symplectic four-manifold,
 symplectic critical surfaces in Kahler surface, respectively. It is because that the Riemann surface immersed in Hermitian surface still has 
the Euler characteristic of tangent bundle, the Euler characteristic of normal bundle and the first Chern class, one can ask: \emph{Do the 
Webster's formulae hold for certain class of closed Riemann surfaces immersed in Hermitian surface}? To give an answer for this problem, 
C.K. Peng and the author (see \cite{[PX-21]}) introduce the Chern-minimal surface and prove that Webster's formulae (\ref{Z-1}) and (\ref{Z-2}) hold
 for the closed generic Chern-minimal surfaces in compact Hermitian surface. In this paper we will study the compactness
of Chern-minimal surfaces.

We first view the Chern-minimal surface from the aspect of analysis. More generally, for a smooth immersion $f$ from Riemann surface $(\Sigma, \emph{\texttt{j}}, ds^2_\Sigma)$ into
Hermitian surface $(M, J, h)$, the tangent map $df$ is a smooth section of bundle $T^*\Sigma\otimes f^{-1}TM$. Then, by using the connection $\nabla$ on $T^*\Sigma\otimes f^{-1}TM$ induced from
Chern connection and taking the covariant differentiation, we obtain the second fundamental form $\nabla df$. We call $f$ is \emph{Chern-harmonic} if and only if it satifies the
equation 
\begin{equation}\label{Z-3}
tr\nabla df=0,
\end{equation}
which is a new elliptic equation involving the metrics and complex structures on manifolds.
Here the trace is taken with respect to the metric $ds^2_\Sigma$. In particular, we call $f$ is \emph{Chern-minimal} if it is an isometry. 
The definitions are analogous to the harmonic map and minimal surface from the viewpoint of geometrical definition. 
It should be pointed out that the Chern-harmonic (resp. Chern-minimal) is just the harmonic (resp. minimal) when the target manifold is K\"ahlerian, and the holomorphic/anti-holomorphic maps
are automatically Chern-harmonic. Although we have not find the variational structure of equation (\ref{Z-3}), it is a conformally invariant equation (see Proposition \ref{Proposition 2.1})
from two-dimensional manifolds. So, it is possible to obtain a compactness result for Chern-minimal surfaces.

There are two celebrated works on compactness of two-dimensional geometric objects, the Sacks-Uhlenbeck's work (see \cite{[SU-81]}) on harmonic maps and M. Gromov's work on pseudo-holomorphic curves.
M. Gromov's original proof is entirely geometric. Inspiring from the works on harmonic maps, T.H. Parker and J.G. Wolfson (see \cite{[PW-93]}), R.G. Ye (see \cite{[Y-94]}) 
independently gave the analytic proof of Gromov's compactness theorem. Their proof reveal the unity of such two compactness problems. Namely, the entire bubble tree procedure requires only a 
conformally invariant elliptic equation with the properties of energy estimate, energy gap and removable sigularity.

Notice that the Chern-harmonic equation is analogous to harmonic map equation and the holomorphic curve is a special case of Chern-harmonic map, so we can study the compactness of Chern-minimal surfaces follows 
the way of harmonic maps. We first deduce a Bochner formula for the energy density of Chern-harmonic map to get the energy estimate, then we prove the isoperimetric inequality and removable sigularity 
theorem for Chern-minimal maps. Based on these preliminaries, we obtain the $C^\infty$-convergence and bubbling.

\vspace{0.2cm}
\noindent\textbf{Theorem 4.1.} \emph{Let $\{f_n\}$ be a sequence of conformal Chern-minimal immersions from closed Riemann surface $(\Sigma, \texttt{j})$ into compact Hermitian surface $(M, J, h)$ 
with the areas are uniformly bounded by $\mathcal{A}_0$. 
Then there are a subsequence $\{f_{n,k}\}$, a finite set of points $\{p_1,\ldots,p_{\kappa}\}\subset\Sigma$ and a conformal Chern-minimal immersion $f_0$ from $\Sigma$ into $M $ such that 
$f_{n,k}$ converges to $f_0$ in $C^\infty$-topology on $\Sigma\setminus \{p_1, \ldots, p_{\kappa}\}$. Furthermore, there is a Chern-minimal immersion $f_{p_i}$ from $S^2$ into $M$ associated to each $p_i$}.

\vspace{0.2cm}
\noindent The proof of Theorem \ref{Theorem 4.1} is essentially due to Sacks and Uhlenbeck (see \cite{[SU-81]}), it is now  known as Sacks-Uhlenbeck's procedure. This procedure shows the existence of Chern-harmonic
surface, but it doesn't characterize the quantities of energy and the position of bubbles. By modifying the renormalization procedure of T.H. Parker and J.G. Wolfson (see \cite{[PW-93]}), T.H. Parker (see \cite{[P-96]}),
working in conformal coordinate instead of the normal coordinate therein, we prove the energy identity and the necklessness, and hence we get the bubble tree convergence. The no energy loss statement for harmonic
maps was previously proved by J. Jost in \cite{[J-91]}.

\vspace{0.2cm}
\noindent\textbf{Theorem 5.3.} \emph{Let $\{f_n\}$ be a sequence of conformal Chern-minimal immersions from Riemann surface
$(\Sigma, \texttt{j})$ into compact Hermitian surface $(M, J, h)$ with the areas are uniformly bounded by $\mathcal{A}_0$.
Then there are a subsequence $\{f_{n,k}\}$, a Chern-minimal immersion $f_0:\Sigma \longrightarrow M$, a finite set of renormalized 
Chern-minimal sequences $\{\tilde{f}_{n,k,I}\}$ and
a finite set of Chern-minimal two-spheres $f_{p_I}:S^2\longrightarrow M$ so that}

\emph{$(1)$ The sequences $\{f_{n,k}\}$, $\{\tilde{f}_{n,k,I}\}$ converge to $f_0$, $f_{p_I}$ in $C^\infty$-topology on $\Sigma\setminus\{p_1,\ldots,p_\kappa\}$, 
$S^2\setminus\{p_{I1},\ldots, p_{I\kappa_I}, \textbf{s}\}$, respectively, where $\textbf{s}$ is the south pole of $S^2$.}

\vspace{0.1cm}

\emph{$(2)$ There is no energy loss. That is}
\begin{equation}\label{D-34}
\lim\limits_{k\rightarrow\infty}E(f_{n,k})=E(f_0)+\sum\limits_{I}E(f_{p_I}).
\end{equation}

\emph{$(3)$ There is no distance bubbling. Namely, for each bubble point $p_I$, we have $f_{p_I}(\textbf{s})=f_{p_{I'}}(p_I)$ 
with indices $I'=i_1\cdots i_{\ell-1}$ and $I=i_1\cdots i_{\ell-1}i_\ell$}.

\vspace{0.2cm}
The paper is organized as follows. In Section 2, we deduce a Bochner formula of Chern-harmonic maps, and use it to get the energy estimate.  
In Section 3, we prove the removable singularity of Chern-minimal surfaces by using the isoperimetric inequality. 
In Section 4, we use the Sacks-Uhlenbeck's procedure to get the $C^\infty$-convergence and bubbling. In Section 5, we prove the energy identity and
necklessness to get the bubble tree limit.

\vspace{0.3cm}

\section{A Bochner formula and energy estimate}
The purpose of this section is to deduce a Bochner formula for the energy density of Chern-harmonic map, and we use it to give the energy estimate of Chern-harmonic maps.

Let $(\Sigma, \emph{\texttt{j}})$ be a Riemann surface, and let $ds^2_\Sigma$ be a Riemannian metric on $\Sigma$ which is conformal to the complex structure 
$\emph{\texttt{j}}$. Locally, we choose a complex valued 1-form $\varphi$ of $(1,0)$-type such that $ds^2_\Sigma=\varphi\;\bar{\varphi}$.
Then the structure
equations of $ds^2_\Sigma$ are given by
\begin{equation}\label{A-1}
d\varphi=-\sqrt{-1}\rho\wedge \varphi,\hspace{0.3cm}\bar{\rho}=\rho,
\end{equation} 
\begin{equation}\label{A-2}
d\rho=-\sqrt{-1}K \varphi\wedge \bar{\varphi},
\end{equation}
where $\rho$ is the connection 1-form and $K$ is the Gaussian curvature.

Let $(M,J,h)$ be a Hermitian surface with the complex structure $J$ and the Hermitian metric $h$. We choose the Chern connection on the tangent bundle $TM$. Namely, 
the unique one preserves $J$, $h$ and it also has vanishing $(1,1)$-part of the torsion. 
Choosing a local field of unitary frame $\{e_1,e_2\}$ with the dual $\{\omega^1,\omega^2\}$ on $M$. Then the structure equations of the Chern 
connection are given by 
 \begin{equation}\label{A-3}
d\omega^i=-\omega^i_j\wedge\omega^j+\Theta^i, \hspace{0.5cm}\omega^j_i+\overline{\omega^i_j}=0,
\end{equation}
\begin{equation}\label{A-4}
d\omega^i_j=-\omega^i_k\wedge\omega^k_j+\Omega^i_j,\hspace{0.5cm}\Omega^j_i+\overline{\Omega^i_j}=0,
\end{equation}
where $\omega^i_j$, $\Theta^i$, $\Omega^i_j$ are connection 1-forms, torsion 2-forms and curvature 2-forms, respectively. Explicitly, we can write
\begin{equation}\label{A-5}
\Theta^i:=L^i_{jk}\,\omega^j\wedge\omega^k, \hspace{0.3cm}L^i_{jk}=-L^i_{kj},
\end{equation}
\begin{equation}\label{A-6}
\Omega^i_j:=R^i_{jk\ell}\,\omega^k\wedge\omega^\ell+R^i_{jk\bar{\ell}}\,\omega^k\wedge\overline{\omega^\ell}+R^i_{j\bar{k}\,\bar{\ell}}\,\overline{\omega^k}\wedge\overline{\omega^\ell},
\end{equation}
with $R^i_{jk\ell}=-R^i_{j\ell k}$, $R^i_{j\bar{k}\,\bar{\ell}}=-R^i_{j\bar{\ell}\,\bar{k}}$, $R^i_{jk\ell}=\overline{R^j_{i\bar{\ell}\,\bar{k}}}$ and $R^i_{jk\bar{\ell}}=\overline{R^j_{i\ell\bar{k}}}$.

Let $f$ be a smooth immersion from surface $(\Sigma, \emph{\texttt{j}})$ into Hermitian surface $(M,J,h)$. 
Set
\begin{equation}\label{A-7}
f^*\omega^i:=a^i_1\,\varphi+a^i_{\bar{1}}\,\bar{\varphi}.
\end{equation}
Taking the exterior differentiation of (\ref{A-7}) and defining $a^i_{11}$, $a^i_{1\bar{1}}$, $a^i_{\bar{1}1}$, $a^i_{\bar{1}\,\bar{1}}$ as follows 
\begin{equation}\label{A-8}
a^i_{11}\,\varphi+a^i_{1\bar{1}}\,\bar{\varphi}:=da^i_1-\sqrt{-1}\rho\, a^i_1 +a^j_1\,\omega^i_j,
\end{equation}
\begin{equation}\label{A-9}
a^i_{\bar{1}1}\,\varphi+a^i_{\bar{1}\,\bar{1}}\,\bar{\varphi}:=da^i_{\bar{1}}+\sqrt{-1}\rho\, a^i_{\bar{1}} +a^j_{\bar{1}}\,\omega^i_j.
\end{equation}
Then, we have
\begin{equation}\label{A-10}
(a^i_{11}\,\varphi+a^i_{1\bar{1}}\,\bar{\varphi})\wedge\varphi+(a^i_{\bar{1}1}\,\varphi+a^i_{\bar{1}\,\bar{1}}\,\bar{\varphi})\wedge\bar{\varphi}=f^*\Theta^i,
\end{equation}
which implies
\begin{equation}\label{A-11}
-a^i_{1\bar{1}}+a^i_{\bar{1}1}=2L^i_{jk}a^j_1a^k_{\bar{1}}.
\end{equation}
We take the exterior differentiation of (\ref{A-8}) and (\ref{A-9}), respectively, and defining $a^i_{111}$, $a^i_{11\bar{1}}$, $a^i_{1\bar{1}1}$,
 $a^i_{1\bar{1}\,\bar{1}}$, $a^i_{\bar{1}11}$, $a^i_{\bar{1}1\bar{1}}$, $a^i_{\bar{1}\,\bar{1}1}$, $a^i_{\bar{1}\,\bar{1}\,\bar{1}}$ as follows
\begin{equation}\label{A-12}
a^i_{111}\varphi+a^i_{11\bar{1}}\,\bar{\varphi}:=da^i_{11}-2a^i_{11}\sqrt{-1}\rho+a^j_{11}\omega^i_j,
\end{equation}
\begin{equation}\label{A-13}
a^i_{1\bar{1}1}\varphi+a^i_{1\bar{1}\,\bar{1}}\,\bar{\varphi}:=da^i_{1\bar{1}}+a^j_{1\bar{1}}\,\omega^i_j,
\end{equation}
\begin{equation}\label{A-14}
a^i_{\bar{1}11}\varphi+a^i_{\bar{1}1\bar{1}}\,\bar{\varphi}:=da^i_{\bar{1}1}+a^j_{\bar{1}1}\omega^i_j,
\end{equation}
\begin{equation}\label{A-15}
a^i_{\bar{1}\,\bar{1}1}\varphi+a^i_{\bar{1}\,\bar{1}\,\bar{1}}\,\bar{\varphi}:=da^i_{\bar{1}\,\bar{1}}+2a^i_{\bar{1}\,\bar{1}}\sqrt{-1}\rho+a^j_{\bar{1}\,\bar{1}}\omega^i_j.
\end{equation}
Then, we have
\begin{equation}\label{A-16}
(a^i_{111}\varphi+a^i_{11\bar{1}}\,\bar{\varphi})\wedge\varphi+(a^i_{1\bar{1}1}\varphi+a^i_{1\bar{1}\,\bar{1}}\,\bar{\varphi})\wedge\bar{\varphi}=-Ka^i_1\varphi\wedge\bar{\varphi}+a^j_1\Omega_j^i,
\end{equation}
\begin{equation}\label{A-17}
(a^i_{\bar{1}11}\varphi+a^i_{\bar{1}1\bar{1}}\,\bar{\varphi})\wedge\varphi+(a^i_{\bar{1}\,\bar{1}1}\varphi+a^i_{\bar{1}\,\bar{1}\,\bar{1}}\,\bar{\varphi})\wedge\bar{\varphi}
=Ka^i_{\bar{1}}\,\varphi\wedge\bar{\varphi}+a^j_{\bar{1}}\,\Omega_j^i,\hspace{0.1cm}{}
\end{equation}
which imply the Ricci identities
\begin{equation}\label{A-18}
a^i_{1\bar{1}1}-a^i_{11\bar{1}}=-Ka^i_1+a^j_1\;\Omega^i_j/(\varphi\wedge\bar{\varphi}),
\end{equation}
\begin{equation}\label{A-19}
a^i_{\bar{1}\,\bar{1}1}-a^i_{\bar{1}1\bar{1}}=Ka^i_{\bar{1}}+a^j_{\bar{1}}\;\Omega^i_j/(\varphi\wedge\bar{\varphi}),
\end{equation}
where 
\begin{equation}\label{A-20}
\Omega^i_j/(\varphi\wedge\bar{\varphi}):=2R^i_{jk\ell}a^k_1a^\ell_{\bar{1}}+R^i_{jk\bar{\ell}}(a^k_1\overline{a^\ell_1}-a^k_{\bar{1}}\overline{a^\ell_{\bar{1}}})
+2R^i_{j\bar{k}\,\bar{\ell}}\overline{a^k_{\bar{1}}}\,\overline{a^\ell_1}
\end{equation}
stands for the coefficient of the pull-back of $\Omega^i_j$ with respect to the 2-form $\varphi\wedge\bar{\varphi}$.

By using (\ref{A-7}), (\ref{A-8}) and (\ref{A-9}), we can give a local expression of the second fundamental form as
\begin{eqnarray}\label{A-21}
\nabla df&=&a^i_{11}\,\varphi\otimes\varphi\otimes e_i+a^i_{1\bar{1}}\,\bar{\varphi}\otimes\varphi\otimes e_i
+a^i_{\bar{1}1}\,\varphi\otimes\bar{\varphi}\otimes e _i+a^i_{\bar{1}\,\bar{1}}\,\bar{\varphi}\otimes\bar{\varphi}\otimes e_i\nonumber\\
&{}&+\,\overline{a^i_{11}}\,\bar{\varphi}\otimes\bar{\varphi}\otimes \overline{e_i}+\overline{a^i_{1\bar{1}}}\,\varphi\otimes\bar{\varphi}\otimes\overline{e_i}
+\overline{a^i_{\bar{1}1}}\,\bar{\varphi}\otimes\varphi\otimes \overline{e_i}+\overline{a^i_{\bar{1}\,\bar{1}}}\,\varphi\otimes\varphi\otimes \overline{e_i}.\hspace{0.6cm}
\end{eqnarray}
So, $f$ is Chern-harmonic if and only if 
\begin{equation}\label{A-22}
a^i_{1\bar{1}}+a^i_{\bar{1}1}=0.
\end{equation}
This together with (\ref{A-11}) imply
\begin{equation}\label{A-23}
a^i_{1\bar{1}}=-a^i_{\bar{1}1}=-L^i_{jk}a^j_1a^k_{\bar{1}}.
\end{equation}

The following proposition shows that Chern-harmonic equation (\ref{Z-3}) is conformally invariant.
\begin{Proposition}\label{Proposition 2.1}
Let $f$ be a smooth map from Riemann surface $(\Sigma, \texttt{j})$ into Hermitian surface $(M, J, h)$. The metrics $ds^2_\Sigma$, $d\tilde{s}^2_\Sigma$ 
are conformal to $\texttt{j}$. Then $f$ is Chern-harmonic with respect to $ds^2_\Sigma$
if and only if it is Chern-harmonic with respect to $ d\tilde{s}^2_\Sigma$.
\end{Proposition}

\emph{Proof}. Locally, we write $h=\omega^1\overline{\omega^1}+\omega^2\overline{\omega^2}$, $ds^2_\Sigma=\varphi\bar{\varphi}$ and $d\tilde{s}^2_\Sigma=\theta\bar{\theta}$, 
and set
\begin{equation}\label{A-24}
f^*\omega^i=a^i_1\varphi+a^i_{\bar{1}}\bar{\varphi},\hspace{1cm} f^*\omega^i=b^i_1\theta+b^i_{\bar{1}}\bar{\theta}.
\end{equation}Notice that both $ds^2_\Sigma$ and $d\tilde{s}^2_\Sigma$ are conformal to $\emph{\texttt{j}}$, there is a local smooth function $\mu$ such that
\begin{equation}\label{A-25}
\varphi=\mu\,\theta.
\end{equation}
We define the convariant derivatives $\mu_1$, $\mu_{\bar{1}}$ by
\begin{equation*}
\mu_1\theta+\mu_{\bar{1}}\,\bar{\theta}:=d\mu-\sqrt{-1}\mu\,\rho_0+\sqrt{-1}\mu\,\rho.
\end{equation*}
Taking the exterior differentiation of (\ref{A-25}), we have
\begin{equation*}
(\mu_1\theta+\mu_{\bar{1}}\,\bar{\theta})\wedge \theta=0,
\end{equation*}
which implies that $\mu_{\bar{1}}=0$. On the other hand, it follows from (\ref{A-24}) and (\ref{A-25}) we have $b^i_1=\mu\,a^i_1$, $b^i_{\bar{1}}=\bar{\mu}\,a^i_{\bar{1}}$.
So, we have
\begin{eqnarray}\label{A-26}
b^i_{1\bar{1}}=(\mu\,a^i_1)_{\bar{1}}&=&\mu_{\bar{1}}\, a^i_{1}+\mu\,(a^i_{1})_{\bar{1}}\nonumber\\
&=&\mu_{\bar{1}}\, a^i_{1}+\mu\,a^i_{1\bar{1}}\,\bar{\mu}\nonumber\\
&=&|\mu|^2\,a^i_{1\bar{1}}.
\end{eqnarray}
Similarly, we also have 
\begin{equation}\label{A-27}
b^i_{\bar{1}1}=|\mu|^2a^i_{\bar{1}1}.
\end{equation}
Then the conclusion follows from (\ref{A-22}), (\ref{A-26}) and (\ref{A-27}).

\hfill$\Box$

We next provide an alternative description to the Chern-harmonic map. Recalling that the differential operator $d^\nabla$ acts on $\phi\otimes e\in \Omega^p(f^{-1}TM)$ 
is defined by $d^\nabla(\phi\otimes e)=d\phi\otimes e+(-1)^p\phi\wedge \nabla e$, 
one can check that $f$ is Chern-harmonic if and only if 
\begin{equation}\label{A-28}
(d^\nabla)^*df=0,
\end{equation}
where $(d^\nabla)^*$ is the adjoint of $d^\nabla$. Moreover, for a Chern-harmonic map $f$, we have 
\begin{equation}\label{A-29}
\big((d^\nabla)^*+d^\nabla\big) df=p(df,df),
\end{equation}
and
\begin{equation}\label{A-30}
\big(d^\nabla(d^\nabla)^*+(d^\nabla)^*d^\nabla\big) df=q(\nabla df,df,df,df),
\end{equation}
where
\begin{equation*}
p(df,df):=2L^i_{jk}a^j_1a^k_{\bar{1}}\;\varphi\wedge\bar{\varphi}\otimes e_i-2\overline{L^i_{jk}}\overline{a^j_1}\overline{a^k_{\bar{1}}}\;\varphi\wedge\bar{\varphi}\otimes\overline{e_i},
\end{equation*}
\begin{eqnarray*}
q(\nabla df,df,df,df)&:=&2\big[(L^i_{jk\ell}a^\ell_1+L^i_{jk\bar{\ell}}\overline{a^\ell_{\bar{1}}})a^j_1a^k_{\bar{1}}+
L^i_{jk}(a^j_{11}a^k_{\bar{1}}+a^j_1a^k_{\bar{1}1})\big]\varphi\otimes e_i\nonumber\\
&{}&-2\big[(L^i_{jk\ell}a^\ell_{\bar{1}}+L^i_{jk\bar{\ell}}\overline{a^\ell_1})a^j_1a^k_{\bar{1}}+L^i_{jk}(a^j_{1\bar{1}}a^k_{\bar{1}}+a^j_1a^k_{\bar{1}\,\bar{1}})\big]\bar{\varphi}\otimes e_i\nonumber\\
&{}&-2\big[(\overline{L^i_{jk\ell}}\,\overline{a^\ell_{\bar{1}}}+\overline{L^i_{jk\bar{\ell}}}a^\ell_1)\overline{a^j_1}\overline{a^k_{\bar{1}}}
+\overline{L^i_{jk}}(\overline{a^j_{1\bar{1}}}\overline{a^k_{\bar{1}}}+\overline{a^j_1}\overline{a^k_{\bar{1}\,\bar{1}}})\big]\varphi\otimes \overline{e_i}\nonumber\\
&{}&+2\big[(\overline{L^i_{jk\ell}}\,\overline{a^\ell_1}+\overline{L^i_{jk\bar{\ell}}}a^\ell_{\bar{1}})\overline{a^j_1}\overline{a^k_{\bar{1}}}+
\overline{L^i_{jk}}(\overline{a^j_{11}}\overline{a^k_{\bar{1}}}+\overline{a^j_1}\overline{a^k_{\bar{1}1}})\big]\bar{\varphi}\otimes \overline{e_i},
\end{eqnarray*}
where $L^i_{jk\ell}$, $L^i_{jk\bar{\ell}}$ are the covariant derivatives of the torsion. It should be pointed out that
$(d^\nabla)^*+d^\nabla$, $d^\nabla(d^\nabla)^*+(d^\nabla)^*d^\nabla$ are elliptic operators of order one and two, respectively.

\begin{Theorem}\label{Theorem 2.2}
Let $f$ be a Chern-harmonic map from Riemann surface $(\Sigma, \texttt{j}, ds^2_\Sigma)$ into Hermitian surface $(M,J,h)$. Then, for the energy density $e(f)$, we have
\begin{eqnarray}\label{A-31}
\Delta e(f)&=&2|A|^2+2K\,e(f)+2\overline{a^i_1}a^j_1\,\Omega^i_j/(\varphi\wedge\bar{\varphi})+2a^i_{\bar{1}}\overline{a^j_{\bar{1}}}\,\Omega^j_i/(\varphi\wedge\bar{\varphi})\nonumber\\
&{}&-4\mbox{Re}\big[a^i_1(\overline{L^i_{jk\ell}}\,\overline{a^\ell_1}+\overline{L^i_{jk\bar{\ell}}}\, a^\ell_{\bar{1}})\overline{a^j_1}\,\overline{a^k_{\bar{1}}}
+a^i_1\overline{L^i_{jk}}(\overline{a^j_{11}}\overline{a^k_{\bar{1}}}+\overline{a^j_1}\overline{a^k_{\bar{1}1}})\big]\nonumber\\
&{}&+4\mbox{Re}\big[a^i_{\bar{1}}(\overline{L^i_{jk\ell}}\,\overline{a^\ell_{\bar{1}}}+\overline{L^i_{jk\bar{\ell}}}a^\ell_1)\overline{a^j_1}\,\overline{a^k_{\bar{1}}}
+a^i_{\bar{1}}\,\overline{L^i_{jk}}(\overline{a^j_{1\bar{1}}}\overline{a^k_{\bar{1}}}+\overline{a^j_1}\overline{a^k_{\bar{1}\,\bar{1}}})\big],
\end{eqnarray}
where $|A|^2=|a^i_{11}|^2+2|a^i_{1\bar{1}}|^2+|a^i_{\bar{1}\,\bar{1}}|^2$ is the half of the squared norm of $\nabla df$.
\end{Theorem}

\emph{Proof}. Notice that the \emph{energy density} $e(f)=|a^i_1|^2+|a^i_{\bar{1}}|^2$, we have
\begin{eqnarray}\label{A-32}
\frac{1}{2}\Delta e(f)&=&(a^i_1\overline{a^i_1}+a^i_{\bar{1}}\overline{a^i_{\bar{1}}})_{1\bar{1}}\nonumber\\
&=&|A|^2+\overline{a^i_1}a^i_{11\bar{1}}+a^i_1\overline{a^i_{1\bar{1}1}}+\overline{a^i_{\bar{1}}}\,a^i_{\bar{1}1\bar{1}}+a^i_{\bar{1}}\,\overline{a^i_{\bar{1}\,\bar{1}1}}.
\end{eqnarray}
It follows from (\ref{A-23}) that
\begin{eqnarray}\label{A-33}
a^i_1\overline{a^i_{1\bar{1}1}}&=&a^i_1(-\overline{L^i_{jk}a^j_1a^k_{\bar{1}}})_{\bar{1}}\nonumber\\
&=&-a^i_1(\overline{L^i_{jk\ell}}\,\overline{a^\ell_1}+\overline{L^i_{jk\bar{\ell}}}\, a^\ell_{\bar{1}})\overline{a^j_1}\,\overline{a^k_{\bar{1}}}
-a^i_1\overline{L^i_{jk}}(\overline{a^j_{11}}\overline{a^k_{\bar{1}}}+\overline{a^j_1}\overline{a^k_{\bar{1}1}}),
\end{eqnarray}
and
\begin{eqnarray}\label{A-34}
\overline{a^i_{\bar{1}}} a^i_{\bar{1}1\bar{1}}=\overline{a^i_{\bar{1}}}(L^i_{jk\ell}a^\ell_{\bar{1}}+L^i_{jk\bar{\ell}}\overline{a^\ell_1})a^j_1a^k_{\bar{1}}
+\overline{a^i_{\bar{1}}}L^i_{jk}(a^j_{1\bar{1}}a^k_{\bar{1}}+a^j_1a^k_{\bar{1}\,\bar{1}}).
\end{eqnarray}
For the other terms in (\ref{A-32}), by using the Ricci identities (\ref{A-18}), (\ref{A-19}) and (\ref{A-23}) again, we have
\begin{eqnarray}\label{A-35}
\overline{a^i_1}a^i_{11\bar{1}}&=&K|a^i_1|^2+\overline{a^i_1}a^j_1\;\Omega^i_j/(\varphi\wedge\bar{\varphi})\nonumber\\
&{}&-\overline{a^i_1}(L^i_{jk\ell}a^\ell_1+L^i_{jk\bar{\ell}}\overline{a^\ell_{\bar{1}}})a^j_1a^k_{\bar{1}}-\overline{a^i_1}L^i_{jk}(a^j_{11}a^k_{\bar{1}}+a^j_1a^k_{\bar{1}1}),
\end{eqnarray}
and 
\begin{eqnarray}\label{A-36}
a^i_{\bar{1}}\overline{a^i_{\bar{1}\,\bar{1}1}}&=&K|a^i_{\bar{1}}|^2+a^i_{\bar{1}}\overline{a^j_{\bar{1}}}\;\overline{\Omega^i_j/(\varphi\wedge\bar{\varphi})}\nonumber\\
&{}&+a^i_{\bar{1}}(\overline{L^i_{jk\ell}}\,\overline{a^\ell_{\bar{1}}}+\overline{L^i_{jk\bar{\ell}}}a^\ell_1)\overline{a^j_1}\,\overline{a^k_{\bar{1}}}
+a^i_{\bar{1}}\,\overline{L^i_{jk}}(\overline{a^j_{1\bar{1}}}\overline{a^k_{\bar{1}}}+\overline{a^j_1}\overline{a^k}_{\bar{1}\,\bar{1}}).
\end{eqnarray}
Substituting (\ref{A-33})-(\ref{A-36}) into (\ref{A-32}), we get the Bochner formula (\ref{A-31}).

\hfill$\Box$

\emph{Remark}. The \emph{energy} $E(f)$ of a smooth map $f$ from Riemann surface $(\Sigma, \emph{\texttt{j}}, ds^2_\Sigma)$ is defined by the integration of its energy density. Namely,
\begin{equation}\label{A-37}
E(f)=\int_{\Omega}e(f)\,dA,
\end{equation}
where $dA$ is the area element of $ds^2_\Sigma$. It is clear that $E(f)$ is equal to the area $\mathcal{A}(f(\Sigma))$ when $f$ is conformal.

\begin{Corollary}\label{Corollary 2.3}
Let $f$ be a Chern-harmonic map from closed Riemann surface $(\Sigma,\texttt{j},ds^2_\Sigma)$ into compact Hermitian surface $(M,J,h)$. Then there are positive constants $C_1$ and $C_2$ so that
\begin{equation}\label{A-38}
\Delta e(f)\geq -C_1\,e(f)-C_2\,e^2(f),
\end{equation}
where $C_1$ depends on the curvature of $ds^2_\Sigma$, $C_2$ depends on the torsion and curvature of $h$. 
\end{Corollary}

\emph{Proof}. It follows from (\ref{A-20}), (\ref{A-31}) and the Cauchy's inequality with $\epsilon$.
\hfill$\Box$

\vspace{0.2cm}

\emph{Remark}. If we scaling the metric $ds^2_\Sigma$ as
$d\tilde{s}^2_\Sigma=\lambda^2 ds^2_\Sigma$ for a positive constant $\lambda$, then the corresponding constants $\tilde{C}_1=C_1/\lambda^2$ and $\tilde{C}_2=C_2$.

\vspace{0.2cm}

Once we have the differential inequality (\ref{A-38}), one can get the following energy estimate of Chern-harmonic maps.
This approach is the same as the corresponding one in the theory of harmonic maps and pseudo-holomorphic curves, which has been
used by R. Schoen \cite{[S-84]}, J.G. Wolfson \cite{[W-88]} and T.H. Parker, J.G. Wolfson \cite{[PW-93]}.

\begin{Theorem}\label{Theorem 2.4}
Let $f$ be a Chern-harmonic map from closed Riemann surface $(\Sigma,\texttt{j},ds^2_\Sigma)$ into compact Hermitian surface $(M,J,h)$. Then there are constants $C_3$, $\epsilon_1>0$  
depend on the metrics $ds^2_\Sigma$ and $h$ so that for any geodesic disk $D_{2r}$ of radius $2r$ with the energy $E(2r):=\int_{D_{2r}}e(f)\,dA\leq \epsilon_1$, we have
\begin{equation}\label{A-39}
\sup\limits_{D_{r}}\,e(f)\leq C_3\,\frac{E(2r)}{r^2}.
\end{equation}
\end{Theorem}

\emph{Proof}. We consider the function $u(\tau):=\tau^2\sup\limits_{D_{2(r-\tau)}}e(f)$ with $\tau_0$ is a maximum point of $u(\tau)$.
Set $e_0:=\sup\limits_{D_{2(r-\tau_0)}}e(f)$, and let $x_0$ be a point in $D_{2(r-\tau_0)}$ such that $e(f)(x_0)=e_0$. It follows from $u(\tau_0/2)\leq u(\tau_0)$, we have
\begin{equation}\label{A-40}
\sup\limits_{D_{2r-\tau_0}} e(f)\leq 4e_0.
\end{equation}
Notice that $D_{\tau_0}(x_0)\subset D_{2r-\tau_0}$, then from (\ref{A-40}) we have $e(f)\leq 4e_0$ on $D_{\tau_0}(x_0)$. Set $d\tilde{s}^2_{\Sigma}=4e_0\,ds^2_\Sigma$, considering the Chern-harmonic map $f$
from $(\Sigma,\emph{\texttt{j}}, d\tilde{s}^2_\Sigma)$ into Hermitian surface $(M,J,h)$, we have $\tilde{e}(f)=(4e_0)^{-1}\,e(f)\leq 1$ on $D_{\tau_0}(x_0)$. By using the fact $\tilde{\Delta}=(4e_0)^{-1}\Delta$,
then the Corollary \ref{Corollary 2.3} and its remark imply
\begin{equation}\label{A-41}
\Delta \tilde{e}(f)\geq -C_1\,\tilde{e}(f)-4C_2e_0\,\tilde{e}^2(f),
\end{equation}
where $C_1$, $C_2$ are the constants in Corollary \ref{Corollary 2.3}.

We first consider the case that $e_0\geq 1$. On $D_{\tau_0}(x_0)$, by the fact $\tilde{e}(f)\leq 1$, then (\ref{A-41}) gives 
\begin{equation}\label{A-42}
\Delta \tilde{e}(f)\geq -C_{3.1}\,e_0, 
\end{equation}
where $C_{3.1}=C_1+4C_2$. Applying the Theorem 9.20 in \cite{[GT-83]} to (\ref{A-42}) on $D_\tau(x_0)$ with $0<\tau\leq\tau_0$, we have
\begin{eqnarray*}
\frac{1}{4}\leq\sup\limits_{D_{\tau/2}(x_0)}\,\tilde{e}(f)&\leq& C_{3.2}\big[\frac{1}{4e_0\tau^2}\int_{D_\tau(x_0)}e(f)\;dA +C_{3.1}\;e_0\tau^2\big]\nonumber\\
&=&\frac{C_{3.2}}{4e_0\tau^2}\,E(2r)+C_{3.1}C_{3.2}\,\tau^2\,e_0,
\end{eqnarray*}
which implies
\begin{equation}\label{A-43}
e_0\leq 4C_{3.1}C_{3.2}\,\tau^2e^2_0+\frac{C_{3.2}}{\tau^2}\,E(2r),
\end{equation}
where $C_{3.2}$ depends on the metric $ds^2_\Sigma$. We claim that $4C_{3.1}C_{3.2}\tau^2_0 e_0<1/2$ for $\epsilon_1<\frac{1}{16C_{3.1}C_{3.2}^2}$. 
If not, suppose that $4C_{3.1}C_{3.2}\tau^2_0 e_0\geq 1/2$ and taking
$\tau=\frac{1}{\sqrt{8C_{3.1}C_{3.2}e_0}}\leq \tau_0$ in (\ref{A-43}), we obtain
\begin{equation*}
e_0\leq \frac{e_0}{2}+8C_{3.1}C_{3.2}^2\,E(2r)\,e_0,
\end{equation*}
which is a contradiction when $E(2r)\leq \epsilon_1$. Taking $\tau=\tau_0$ in (\ref{A-43}) and using $4C_{3.1}C_{3.2}\tau^2_0 e_0<1/2$ for $\epsilon_1<\frac{1}{16C_{3.1}C_{3.2}^2}$, we have
\begin{equation}\label{A-44}
u(\tau_0)=\tau_0^2\sup\limits_{D_{2(r-\tau_0)}}e(f)=\tau_0^2e_0\leq 2C_{3.2}\,E(2r).
\end{equation}
So, (\ref{A-39}) follows from $u(r/2)\leq u(\tau_0)$ with $C_3=8C_{3.2}$.

We next consider the case that $e_0<1$. By the fact $\tilde{e}(f)\leq 1$ on $D_{\tau_0}(x_0)$, (\ref{A-41}) gives
\begin{equation}\label{A-45}
\Delta \tilde{e}(f)\geq -C_{3.1}\,\tilde{e}(f),
\end{equation}
where $C_{3.1}=C_1+4C_2$ as before. Applying the Theorem 9.20 in \cite{[GT-83]} again, (\ref{A-45}) gives
\begin{eqnarray*}
\frac{1}{4}\leq\sup\limits_{D_{\tau/2}(x_0)}\,\tilde{e}(f)&\leq& \frac{C_{3.4}}{4\tau^2_0e_0}\int_{D_\tau(x_0)}e(f) \;dA \leq  \frac{C_{3.4}}{4\tau^2_0e_0}\,E(2r),
\end{eqnarray*}
which implies 
\begin{equation}\label{A-46}
u(\tau_0)=\tau_0^2\sup\limits_{D_{2(r-\tau_0)}}e(f)=\tau_0^2e_0\leq C_{3.4}\,E(2r).
\end{equation}
Then the conclusion follows from $u(r/2)\leq u(\tau_0)$ with $C_3=4C_{3.4}$.

\hfill$\Box$
\vspace{0.3cm}

By the definition of Chern-minimal, we know that $ds^2_\Sigma$ is equal to the pull-back metric $f^*h$, which means that $ds^2_\Sigma$ depends on $f$. 
However, one need the uniformly estimates when considering the compactness problem of conformal Chern-minimal maps. 
So, it is necessary to choose a good background metric on $(\Sigma,\emph{\texttt{j}})$.
According to the uniformization theorem, we can always choose the metric $ds_0^2$ on $\Sigma$ so that it
is conformal to $\emph{\texttt{j}}$ and it has constant curvature $1$, $0$, $-1$ when the genus $g(\Sigma)=0$, $1$, $g(\Sigma)\geq 2$, respectively.

\vspace{0.3cm}

\section{Isoperimetric inequality and removable singularity}
The goal of this section is to prove the removable sigularity theorem for Chern-minimal maps, which is based on an isoperimetric inequality, Morrey's decay lemma, energy estimate and elliptic estimates.

\begin{Theorem}\label{Theorem 3.1}
 Let $f$ be a Chern-minimal immersion from Riemann surface $(\Sigma, \texttt{j}, ds^2_\Sigma)$ into compact Hermitian surface $(M, J, h)$. Then there are a universal positive constant $C_4$
and a positive constant $\epsilon_2$ depends on $h$ such that for any domain $\Omega\subset \Sigma$ with boundary provided that $\mathcal{A}(f(\Omega))\leq\epsilon_2$, we have 
\begin{equation}\label{B-1}
\mathcal{A}(f(\Omega))\leq C_4\,\mathcal{L}^2(f(\partial \Omega)),
\end{equation}
where $\mathcal{L}(f(\partial \Omega))$ is the length of $f(\partial \Omega)$.
\end{Theorem}

\emph{Proof}. We denote by $H$ the mean curvature of $f$ in the sense of Levi-Civita connection of $ds^2_\Sigma$ and $h$. By the Proposition 2.2 in \cite{[PX-21]}, we have
\begin{equation}\label{B-2}
H=2(a^j_{\bar{1}}\overline{L^j_{ki}}\overline{a^k_{\bar{1}}}+\overline{a^j_1}\overline{L^k_{ji}}a^k_1)e_i
+2(\overline{a^j_{\bar{1}}}L^j_{ki}a^k_{\bar{1}}+a^j_1L^k_{ji}\overline{a^k_1})\overline{e_i}.
\end{equation}
Notice that $f$ is an isometry, we know that $|a^i_1|, |a^i_{\bar{1}}|\leq 1$. So, from (\ref{B-2}), we have $|H|\leq C_{4.1}$ for a positive constant $C_{4.1}$ depends on the 
torsion of Chern connection of $h$. Then, it follows from the Theorem 2.2 in \cite{[HS-74]} that
\begin{equation}\label{B-3}
\mathcal{A}^{1/2}(f(\Omega))\leq C_{4.2} \Big(\mathcal{L}(f(\partial \Omega))+\int_{\Omega}|H|\,dA\Big)
\end{equation}
provided that $\mathcal{A}(f(\Omega))\leq C_{4.3}$, where $C_{4.2}$ is a universal positive constant and $C_{4.3}$ is a positive constant depends on the injectivity radius and sectional curvature of $h$. 
Then, the conclusion (\ref{B-1}) follows from (\ref{B-3}) with $2C_{4.1}C_{4.2}\mathcal{A}(f(\Omega))\leq \mathcal{A}^{1/2}(f(\Omega))$. 
Namely, we can choose $\epsilon_2<\min\{C_{4.3},\frac{1}{4C_{4.1}^2C_{4.2}^2}\}$ and $C_4=4C^2_{4.2}$.

\hfill$\Box$

\emph{Remark}. Essentially, the isoperimetric inequality holds for Chern-minimal maps is because of the fact that a Chern-minimal surface in compact Hermitian surface has bounded mean curvature.

\begin{Corollary}\label{Corollary 3.2}
 Let $f$ be a Chern-minimal immersion from Riemann suface $(\Sigma, \texttt{j})$ into compact Hermitian surface $(M, J, h)$. Then, for any $x\in f(\Sigma)$ and sufficient small ball $B_r(x)\subseteq M$ with
no boundary points in $B_r(x)$, we have 
\begin{equation}\label{B-4}
C_5\,r^2\leq\mathcal{A}(f(\Sigma)\cap B_r(x)),
\end{equation}
where $C_5=1/4C_4$.
\end{Corollary}

\emph{Proof}. Set $\mathcal{A}(t):=\mathcal{A}(f(\Sigma)\cap B_t(x))$. It follows from Theorem \ref{Theorem 3.1} we obtain $\sqrt{\mathcal{A}(t)}\leq \sqrt{C_4}\, \mathcal{A}'(t)$, whose integration from $0$ to 
$r$ yeilds (\ref{B-4}).

\hfill$\Box$ 

\begin{Corollary}\label{Corollary 3.3}
Let $f$ be a conformal Chern-harmonic map from Riemann surface $(\Sigma, \texttt{j})$ into compact Hermitian surface $(M,J,h)$. Then $f$ is a constant map
if $E(f)\leq \epsilon_2$.
\end{Corollary}

\emph{Proof}. It follows from the remark of Theorem \ref{Theorem 2.2} and Theorem \ref{Theorem 3.1}, for any geodesic disk $D_r$ in $\Sigma$, we have
\begin{equation*}
E(f|_{\Sigma\setminus D_r})=\mathcal{A}(f(\Sigma\setminus D_r))\leq C_4\,\mathcal{L}^2(f(\partial D_r)).
\end{equation*}
Letting $r\rightarrow 0$, we obtain $E(f)=0$, which implies that $f$ is a constant map.

\hfill$\Box$

Inspiring from (\ref{A-24}) and integration by parts, we say a map $f\in W^{1,2}(\Sigma, M)$ is \emph{weakly Chern-minimal} if it satisfies
\begin{equation}\label{B-5}
\int_\Sigma \langle d^\nabla\xi, df\rangle\, dA_0=0,
\end{equation}
for all $\xi\in \Omega^0(f^{-1}TM)$ with compact support. Notice that the regularity problem is a local problem, so we always work on the geodesic disk $D_r$ with the center $p$,
and the punctured disk will be denoted by $D_r^*$.  

\begin{Lemma}\label{Lemma 3.4}
Let $f$ be a smooth conformal Chern-minimal map from a punctured geodesic disk $D^*$ into Hermitian surface $(M, J, h)$. Suppose that $f$ is continuous on $D$ with finite area,
then $f$ is a weakly Chern-minimal map from the geodesic disk $D$ into $M$.
\end{Lemma}

\emph{Proof}. It is sufficient to show that the identity (\ref{B-5}) holds for any $\xi \in \Omega^0(f^{-1}TM)$ with $\mbox{supp}(\xi)\subset D$. For any $0<\epsilon<1$, we 
take a cut-off function $\eta_\epsilon$ such that
\begin{equation*}
\mbox{supp}(\eta_\epsilon)\subset D_\epsilon,\hspace{0.5cm} 0\leq \eta_\epsilon\leq1,\hspace{0.5cm} \eta_\epsilon|_{D_{\epsilon/2}}=1,\hspace{0.5cm} |d\eta_\epsilon|\leq \frac{C'}{\epsilon},
\end{equation*}
where $C'$ is a uniform positive constant. Then
\begin{equation}\label{B-6}
\int_{D}\langle d^\nabla\xi,df\rangle\, dA_0=\int_D\langle d^\nabla((1-\eta_\epsilon)\xi), df\rangle\, dA_0+\int_D\langle d^\nabla(\eta_\epsilon \xi), df\rangle\, dA_0.
\end{equation}
The first term in (\ref{B-6}) is equal to zero by the divergence theorem and (\ref{A-28}). For the last term, notice that 
$d^\nabla(\eta_\epsilon\xi)=d\eta_\epsilon\otimes\xi+\eta_\epsilon \nabla\xi$, then the H\"older inequality and (\ref{A-46}) imply
\begin{equation}\label{B-7}
\int_D\langle d^\nabla(\eta_\epsilon \xi), df\rangle dA_0\leq C'' (\epsilon\,|\nabla\xi|_{L^\infty}+|\xi|_{L^\infty}) \;\mathcal{A}^{1/2}(f(D_\epsilon)),
\end{equation}
where $C''$ is a positive constant depends on $ds^2_0$. So, the last term vanishes by letting $\epsilon\rightarrow 0$ in (\ref{B-7}). This shows that $f$ is a weakly Chern-minimal map from $D$ into $M$.

\hfill$\Box$

\emph{Remark}. The definition of weakly Chern-minimal and the proof of Lemma \ref{Lemma 3.4} are the same as the case of pseudo-holomorphic curves, which are given by T.H. Parker and J.G. Wolfson in \cite{[PW-93]}.

We now prove the \emph{removable singularity theorem} of Chern-minimal maps.

\begin{Theorem}\label{Theorem 3.5}
Let $f$ be a smooth conformal Chern-minimal map from the punctured geodesic disk $D^*$ into compact Hermitian surface $(M, J, h)$ with finite area. 
Then $f$ can be extended to a smooth Chern-minimal map from $D$ into $(M, J, h)$.
\end{Theorem}

\emph{Proof}. We first use the energy estimate and Morrey's decay lemma to prove that $f$ is $C^\alpha$ on $D_{r_0}^*$ for some $\alpha\in (0,1)$ and $r_0\ll 1$, 
and hence $f$ can be extended continuously to $D$, then we use the elliptic estimates
to show that $f$ is smooth on $D$.

Notice that $\mathcal{A}(f(D^*))$ is finite, we have
\begin{equation}\label{B-8}
\lim\limits_{r\rightarrow 0} \mathcal{A}(f(D^*_r))=0.
\end{equation}
Thus, there is a constant $r_1>0$ such that $D^*_{r_1}\subset D^*$ and $\mathcal{A}(f(D^*_{r_1}))=E_0(f|_{D^*_{r_1}})\leq \mbox{min}\{\epsilon_1,\epsilon_2\}$,
where $\epsilon_1, \epsilon_2$ are the constants in Theorem \ref{Theorem 2.4}
and Theorem \ref{Theorem 3.1}, respectively. Let $D_r(x)$ be a geodesic disk with the center $x$ and radius $r$. Obviously, the isoperimetric inequality holds for the domain $D_r(x)\subset D^*_{r_1}$. 
We will show the isoperimetric inequality still holds for the domain $D^*_r(x):=D_r(x)\setminus \{p\}\subset D^*_{r_1}$, where $p$ is the center of $D$.
It follows from Proposition \ref{Proposition 2.1}, we know that $f$ is a Chern-harmonic map from $(D^*, ds^2_0)$ into $(M, J, h)$. So, by the Theorem \ref{Theorem 2.4}, for any $x\in \partial D_\rho$ with $0<2\rho<r_1$, we have
\begin{equation}\label{B-9}
|df|(x)=\sqrt{e_0(f)(x)}\leq \sup\limits_{D_{\rho/2}(x)}\sqrt{e_0(f)}\leq \frac{\sqrt{C_3}}{\rho} \sqrt{\mathcal{A}(f(D^*_{2\rho}))}.
\end{equation}
By using the polar coordinate $(\rho, \theta)$ of $D$ and (\ref{B-9}), we have
\begin{eqnarray}\label{B-10}
\mathcal{L}(f(\partial D_\rho))\leq\int_0^{2\pi}|df(\partial_\theta)|\,d\theta
\leq  C_1'\int_0^{2\pi} \rho\, |df|\,d\theta\leq 2\pi C_1'\sqrt{C_3}\sqrt{\mathcal{A}(f(D^*_{2\rho}))},
\end{eqnarray}
where $C_1'$ depends on $ds^2_0$. We choose $\rho$ sufficient small such that $D^*_\rho\subset D^*_r(x)$. Taking 
$\Omega= D^*_r(x)\setminus\overline{D^*_\rho}$ in Theorem 3.1, we have
\begin{equation}\label{B-11}
\mathcal{A}(f(D^*_r(x)))-\mathcal{A}(f(D^*_\rho))\leq C_4\big[\mathcal{L}\big(f(\partial D_r(x))\big)+\mathcal{L}\big(f(\partial D_\rho)\big)\big]^2.
\end{equation}
Letting $\rho\rightarrow0$ in (\ref{B-11}) and using (\ref{B-8}), (\ref{B-10}), we have
\begin{equation}\label{B-12}
\mathcal{A}(f(D^*_r(x)))\leq C_4\,\mathcal{L}^2\big(f(\partial D_r(x))\big).
\end{equation}

We now fix $r_0$ such that $0<2r_0< r_1$. Let $D_r(x)$ be any geodesic disk is contained in $D^*_{r_0}$, and set $\alpha'=\log r_0/\log(r_0/2)<1$.  
It is easy to check that $D_{r^{\alpha'}}(x)$ or $D^*_{r^{\alpha'}}(x)$ is contained in $D^*_{2r_0}$. Denote by $\mathcal{A}(\rho)$ the area of $D_\rho(x)$ or $D^*_\rho(x)$.  
In the polar coordinate $(\rho, \theta)$ of $D_\rho(x)$ or $D^*_\rho(x)$, the isoperimetric inequality and H\"older inequality imply
\begin{eqnarray}\label{B-13}
\mathcal{A}(\rho)&\leq& C_3\,\Big(\int_0^{2\pi} |df(\partial_\theta)|\;d\theta\Big)^2\nonumber\\
&\leq& C_3\, \Big(\int_0^{2\pi\,}\rho \,w(\rho)\,|df|\;d\theta\Big)^2\nonumber\\
&\leq& 2\pi C_3\,\rho^2 \,w^2(\rho)\,\int_0^{2\pi}\,|df|^2\,d\theta\nonumber\\
&\leq& 2\pi C_3 C_1'\,\rho\,\frac{d}{d\rho}\mathcal{A}(\rho),
\end{eqnarray} 
where $C_1'$ is a positive constant as in (\ref{B-10}). Here $w(\rho)=\sin\rho$, $1$, $\sinh\rho$ if the curvature of background metric $ds^2_0$ is $1$, $0$, $-1$, respectively.  Hence, we have
\begin{equation}\label{B-14}
\frac{(2\pi\,C_3\,C_1')^{-1}}{\rho}\leq\frac{d}{d\rho}\log \mathcal{A}(\rho),
\end{equation}
Integrating (\ref{B-14}) from $r$ to $r^{\alpha'}$ gives
\begin{equation}\label{B-15}
\mathcal{A}(r)=\mathcal{A}\big(f(D_r(x))\big)\leq C'_2\; r^\alpha,
\end{equation}
where $C'_2=e^{(2\pi\,C_3\,C_1')^{-1}}\mathcal{A}(D^*_{2r_0})$ and $\alpha=1-\alpha'$. Thus the Morrey's decay lemma (see Theorem 3.5.2 in \cite{[M-66]}
 or Lemma 2.1.10 in \cite{[LW-08]}) gives that $f$ is $C^\alpha$ on $D^*_{r_0}$, and hence
$f$ can be extended continuously to $D$.

Notice that $f$ is a weakly Chern-minimal by the Lemma \ref{Lemma 3.4} and the local expression of (\ref{B-5}) has the same form as (8.4.1) in \cite{[J-11]}, then 
the Lemma 8.4.3 in \cite{[J-11]} gives $f\in W^{2,2}(D, M)$, and hence $df\in W^{1,2}(D, M)$. The Sobolev embedding implies that $df\in L^q(D, M)$. The index $q$ 
here and below stands for different constants which are strictly greater than 2. So, the term $p(df, df)$ in the right hand side of (\ref{A-29}) belongs to $ L^q(D, M)$ . The elliptic estimates 
implies that $df\in W^{1,q}(D, M)$, and hence 
the term $q(\nabla df, df, df, df)$ in the right hand side of (\ref{A-30}) belongs to $L^q(D, M)$. It follows from (\ref{A-30}) and the elliptic estimates that $df\in W^{3,q}(D, M)$.
Thus, the bootstrapping arguments show that $f\in C^\infty(D, M)$.

\hfill$\Box$

\vspace{0.3cm}

\section{$C^\infty$-convergence and bubbling}

In this section we use the Sacks-Uhlenbeck's procedure (see \cite{[SU-81]}) to get $C^\infty$-convergence of conformal Chern-minimal surface and existence of Chern-minimal two-sphere in Hermitian surface
around the bubble point.

\begin{Theorem}\label{Theorem 4.1}
Let $\{f_n\}$ be a sequence of conformal Chern-minimal immersions from closed Riemann surface $(\Sigma, \texttt{j})$ into compact Hermitian surface $(M, J, h)$ with the areas are uniformly bounded by $\mathcal{A}_0$. 
Then there are a subsequence $\{f_{n,k}\}$, a finite set of points $\{p_1,\ldots,p_{\kappa}\}\subset\Sigma$ and a conformal Chern-minimal immersion $f_0$ from $\Sigma$ into $M $ such that 
$f_{n,k}$ converges to $f_0$ in $C^\infty$-topology on $\Sigma\setminus \{p_1, \ldots, p_{\kappa}\}$. Furthermore, there is a Chern-minimal immersion $f_{p_i}$ from $S^2$ into $M$ associated to each $p_i$.
\end{Theorem}

\emph{Proof}. We choose a constant $r_0>0$ and set $r_m:=2^{-m}r_0$ with $m\in \textbf{N}$. For each $m$, choosing a finite covering $\mathcal{C}_m=\{D_{r_m}(p_\alpha)\}$ of $\Sigma$ such that
each point in $\Sigma$ is covered at most $N$ times by disks in $\mathcal{C}_m$ and $\{D_{r_m/2}(p_\alpha)\}$ is still a covering of $\Sigma$.
It clear that $N$ only depends on $\Sigma $. So, for each $n$, we have
\begin{equation}\label{C-1}
\sum\limits_{\alpha}\int_{D_{r_m}(p_\alpha)} e_0(f_n)\, dA_0\leq N\mathcal{A}_0.
\end{equation}
This implies that there are at most $N\mathcal{A}_0/\epsilon_1$ disks in $\mathcal{C}_m$, on which
\begin{equation}\label{C-2}
\int_{D_{r_m}(p_\alpha)} e_0(f_n)\, dA_0>\epsilon_1,
\end{equation}
where $\epsilon_1$ is the constant in Theorem \ref{Theorem 2.4}. The center points of these disks are at most $N\mathcal{A}_0/\epsilon_1$ sequences of points in $\Sigma$.
Notice that $\mathcal{C}_m$ is a finite covering and $\Sigma$ is compact, we may assume that these center points are fixed by passing to a subsequence of $\{f_n\}$.
For each $m$, we call these center points $\{p_{1,m},\ldots, p_{\kappa,m}\}$ with $\kappa\leq N\mathcal{A}_0/\epsilon_1$. By the Theorem \ref{Theorem 2.4}, the elliptic estimates 
and the Arzela-Ascoli theorem, we can successively choose a subsequence of $\{f_n\}$ that converges in $C^\infty$-topology on every disk $D_{r_m/2}(p_\alpha)$ for
each $D_{r_m}(p_\alpha)\in \mathcal{C}_m\setminus \mathcal{C}_m'$, where $\mathcal{C}_m':=\{D_{r_m}(p_{1,m}),\ldots,D_{r_m}(p_{\kappa,m})\}$. Since $\Sigma$ is compact,
after choosing a subsequence of $\{m\}$, we can assume that $p_{1,m},\ldots, p_{\kappa,m}$ converge to points $p_1,\ldots,p_{\kappa}$ as $m\rightarrow \infty$, respectively.
By a diagonal argument, there is a subsequence of $\{f_n\}$ converges in $C^\infty$-topology on $\Sigma\setminus\{p_1,\ldots,p_{\kappa}\}$, and the limit is denoted
by $f_0$. By the elliptic regularity, $f_0$ is smooth and Chern-minimal on $\Sigma\setminus\{p_1,\ldots,p_{\kappa}\}$. Then the Theorem \ref{Theorem 3.5} implies that $f_0$ can be extended
smoothly to $\Sigma$. 

\vspace{0.1cm} 

We next use the Sacks-Uhlenbeck's procedure to get a bubble for each $p\in\{p_1,\ldots, p_{\kappa}\}$.
For notational convenience, we still use $\{f_n\}$ to relabel the subsequence obtained above.  
Suppose that $r_0>0$ such that $2r_0$ less than the injectivity radius of $(\Sigma, ds^2_0)$ and no other $p_i$ in the geodesic disk $D_{r_0}(p)$. Set
\begin{equation*}
b_n:=\sup\limits_{q\in D_{r_0}(p)}\;|df_n(q)|.
\end{equation*}
It is clear that $\{b_n\}$ is unbounded.
Otherwise, the elliptic estimates and Arzela-Ascoli theorem imply that a subsequence converges on $D_{r_0}(p)$ in $C^\infty$-topology. So, without loss of generality, we can assume that 
$b_n\rightarrow \infty$ as $n\rightarrow \infty$. Let $q_n\in \overline{D_{r_0}(p)}$ be the point such that $|df_n(q_n)|=b_n$. It is easy to check that $q_n\rightarrow p$ as
$n\rightarrow \infty$. 

\vspace{0.1cm}

Locally, we choose a conformal coordinate $(U, \psi; x)$ around $p$ in $D_{r_0}(p)$ such that $\psi(U)=B_2(0)$ with $\psi(p)=0$, where $B_2(0)\subset \mathbb{C}$ 
is the ball centered at $0$ with radius 2. The local expression of metric $ds^2_0$ on $B_2(0)$ can be written as $ds^2_0=\lambda^2(x)\,dxd\bar{x}$. 
For sufficient large $n$, we assume that $\psi(q_n)=x_n$, $B_1(x_n)\subset B_2(0)$ and define
\begin{equation}\label{C-3}
\tilde{f}_n: B_{b_n}(0)\longrightarrow M,\;\; y\mapsto f_n\circ\psi^{-1}(x_n+y/b_n).
\end{equation}
We equip $B_{b_n}(0)$ with the metric $ds^2_n=\lambda_n^2(y)\,dyd\bar{y}$ with $\lambda_n(y)=\lambda(x_n+y/b_n)$. Then $\tilde{f}_n$ is a Chern-harmonic map from $(B_{b_n}(0), ds^2_n)$ 
into $(M, J, h)$ with $|d\tilde{f}_n|_n\leq 1$ on $B_{b_n}(0)$ and $|d\tilde{f}(0)|_n=1$. 
By identifying $S^2\setminus\{\emph{\textbf{s}}\}$ with $\mathbb{C}$ via the southern stereographic projection, where $\emph{\textbf{s}}=(0,0,-1)$ is the south pole, 
we can regards $\tilde{f}_n$ as maps from domains in $S^2\setminus\{ \emph{\textbf{s}} \}$ into $M$. Notice that 
$ds^2_n$ converges smoothly to $ds^2_\infty=\lambda^2(0)\,dyd\bar{y}$ on $\mathbb{C}$ and they are conformal, so Proposition \ref{Proposition 2.1} implies that $\tilde{f}_n$ is Chern-harmonic from 
$(B_{b_n}(0), ds^2_\infty)$ into $(M, J, h)$ with finite area. 
Since $ds^2_n$ and $ds^2_\infty$ are uniformly equivalent, $\{|d\tilde{f}_n|_\infty\}$ is bounded. 
By using a sequence of compact sets exhaust $S^2\setminus\{\emph{\textbf{s}}\}$ and diagonal argument, the elliptic estimates yields a subsequence of $\{\tilde{f}_n\}$ 
that converges in $C^\infty$-topology to a Chern-harmonic map $\tilde{f}$ from $(S^2\setminus\{\emph{\textbf{s}}\}, ds^2_\infty)$ 
into $(M, J, h)$ with finite area. It is also Chern-harmonic with respect to the standard metric on $S^2$ restricted to $S^2\setminus\{\emph{\textbf{s}}\}$. 
The Theorem \ref{Theorem 3.5} gives a Chern-harmonic map $\tilde{f}$ from $S^2$ into $M$. By using the Proposition \ref{Proposition 2.1} again, with respect to the induced metric on $S^2$,
we obtain a Chern-minimal map $f_p$ from $S^2$ into $(M, J, h)$, which is not a constant map by the fact $|d\tilde{f}(0)|_\infty=1$.  

\hfill$\Box$

\emph{Remark}. By passing to a subsequence (will be still relabeled as $\{f_n\}$) in Theorem \ref{Theorem 4.1}, we can define
\begin{equation}\label{C-4}
m_i:=\lim\limits_{r\rightarrow 0}\limsup\limits_{n\rightarrow\infty}\int_{D_{r}(p_i)} e_0(f_n)\,dA_0.
\end{equation}
It follows from Theorem \ref{Theorem 4.1} that we get the measures convergence
\begin{equation}\label{C-5}
e_0(f_n)\,dA_0\rightarrow e_0(f_0)\,dA_0+\sum\limits_{i=1}^{\kappa}m_i\,\delta_{p_i},
\end{equation}
where $e_0(f_n)\, dA_0$ is viewed as measures on $\Sigma$, $\delta_{p_i}$ is the point measure.

\vspace{0.3cm}

\section{The bubble tree convergence}

In this section we give another renormalization procedure to construct bubbles by controlling the energy on each bubble. 
This procedure can be iterated and terminates in finite steps, which is now known as bubble tree convergence in the theory
of harmonic maps and pseudo-holomorphic curves. 

Choosing a conformal coordinate $(U; x)$ around any $p_i\in\{p_1,\ldots,p_{\kappa}\}$ such that there is no other bubble points in $U$ with 
\begin{equation*}
x(p_i)=0,\hspace{0.8cm} x(U)=B_4(0)\subset\mathbb{C},\hspace{0.8cm} ds^2_0=\lambda^2 dx d\bar{x},
\end{equation*}
and we will identify $B_4(0)$ with $U$ in the sequel. Let $\{f_n\}$ be the subsequence has been chosen in Theorem \ref{Theorem 4.1}, and set 
\begin{equation}\label{D-1}
r_n:=\sup \Big\{\;r\;\Big|\;\,\int_{B_{2r}(0)}e_0(f_0)\,dA_0\leq\frac{m_i}{32n^2}\;\mbox{with}\;r\leq\frac{1}{n}\Big\}.
\end{equation}
By (\ref{C-5}) and (\ref{D-1}), we choose a subsequence $\{f_{n,k}\}$ inductively such that
\begin{equation}\label{D-2}
(1-\frac{1}{16k^2})\,m_i\leq\int_{B_{r_k/16k^2}(0)}e_0(f_{n,j})dA_0\leq (1+\frac{1}{16k^2})\,m_i,
\end{equation}
\begin{equation}\label{D-3}
\int_{B_{2r_k}(0)\setminus B_{r_k/16k^2}(0)}e_0(f_{n,j})dA_0\leq\frac{m_i}{16k^2},
\end{equation}
for any $j\geq k$.
By the $C^\infty$-convergence in Theorem \ref{Theorem 4.1} we can further assume that
\begin{equation}\label{D-4}
\sup\limits_{x\in \partial B_{2r_k}(0)}\mbox{dist}\big(f_{n,j}(x),f_0(x)\big)\leq \frac{1}{k},
\end{equation}
\begin{equation}\label{D-5}
\sup\limits_{B_{2r_1}(0)\setminus B_{r_k/16k^2}(0)}\big|e_0(f_{n,j})-e_0(f_0)\big|\leq 1,
\end{equation}
for any $j\geq k$.
We define the center of mass of the measure $e_0(f_{n,k})\,dA_0$ on $B_{2r_k}(0)$ by
\begin{equation}\label{D-6}
\tilde{x}_k:=\int_{B_{2r_k}(0)}x\,e_0(f_{n,k})dA_0\Big/\int_{B_{2r_k(0)}}e_0(f_{n,k})dA_0.
\end{equation}
It follows from (\ref{D-2}), we have
\begin{eqnarray*}
\int_{B_{2r_k}(0)}|x|e_0(f_{n,k}) dA_0&\leq&\frac{r_k}{16k^2}\int_{B_{r_k/16k^2(0)}}e_0(f_{n,k})dA_0+2r_k\int_{B_{2r_k}(0)\setminus B_{r_k/16k^2}(0)}e_0(f_{n,k})dA_0\\
&\leq&(\frac{3}{16k^2}+\frac{1}{256k^4})\,r_km_i,
\end{eqnarray*}
which implies 
\begin{equation}\label{D-7}
|\tilde{x}_k|\leq \frac{r_k}{4k^2}, 
\end{equation}
for sufficient large $k$. Set 
\begin{equation}\label{D-8}
\mu_k:=\sup\Big\{\;\mu\;\Big|\;\,\int_{B_{2r_k}(0)\setminus B_\mu(\tilde{x}_k)}e_0(f_{n,k})dA_0\geq C_0\Big\},
\end{equation}
where $0<C_0< \varepsilon_1/2$ is a fixed renormalization constant which is small enough to be chosen later. It follows from
(\ref{D-7}) we have $B_{r_k/16k^2}(0)\subset B_{r_k/k^2}(\tilde{x}_k)$. So, for sufficient large $k$, (\ref{D-3}) gives
\begin{eqnarray*}
\int_{B_{2r_k}(0)\setminus B_{r_k/k^2}(\tilde{x}_k)}e_0(f_{n,k})dA_0&\leq&\int_{B_{2r_k}(0)\setminus B_{r_k/16k^2}(0)}e_0(f_{n,k})dA_0
< C_0.
\end{eqnarray*}
By the definition of $\mu_k$ in (\ref{D-8}), we have
\begin{equation}\label{D-9}
\mu_k\leq \frac{r_k}{k^2},
\end{equation} 
for sufficient large $k$.

Notice that the Theorem \ref{Theorem 2.4} and a covering argument imply that there is a positive constant $C_1'$ depends on $(\Sigma, ds^2_0)$ so that
\begin{equation}\label{D-10}
\sup\limits_{\Sigma} |df_0|^2\leq C_1'\,E_0(f_0).
\end{equation}
Then (\ref{D-4}) and (\ref{D-10}) give
\begin{eqnarray}\label{D-11}
\mbox{dist}\big(f_{n,k}(\partial B_{2r_k}(0)),f_0(p_i)\big)&\leq& 4r_k\sup\limits_{\Sigma} |df_0|+\sup\limits_{x\in \partial B_{2r_k}(0)}\mbox{dist}\big(f_{n,k}(x),f_0(x)\big)\nonumber\\
&\leq&\frac{C_2'}{k},
\end{eqnarray}
where $C_2'$ depends on $C_1'$ and $E_0(f_0)$. 

For each $k$, we define the conformal transformation on $\mathbb{C}$ as
\begin{equation*}
T_k: \mathbb{C}\longrightarrow \mathbb{C}, \;x\mapsto T_k(x)=\frac{1}{\mu_k}(x-\tilde{x}_k).
\end{equation*}
Let $\pi$ be the south stereographic projection from $S^2\setminus \{\textbf{\emph{s}}\}$ to $\mathbb{C}$, where 
$\textbf{\emph{s}}=(0,0,-1)$ is the south pole. Set $\Omega_k:=\pi^{-1}\circ T_k(B_{2r_k}(0))$, $D_k:=\pi^{-1}\circ T_k(B_{k\mu_k}(\tilde{x}_k))\subset S^2$.
It is clear that $\Omega_k$ and $D_k$ exhaust $S^2\setminus\{\textbf{\emph{s}}\}$ as $k$ goes to infinity, respectively.
Then, define the renormalized maps as
\begin{equation*} 
\tilde{f}_{n,k}:=f_{n,k}\circ T_k^{-1}\circ\pi:\Omega_k\longrightarrow M,
\end{equation*}
whose image are agree with $f_{n,k}|_{B_{2r_k}(0)}$. Notice that $T^{-1}_k$ and $\pi$ are conformal, so $\tilde{f}_{n,k}$ is also Chern-harmonic with respect to 
the standard metric on $S^2$. It follows from (\ref{C-5}), (\ref{D-2}), (\ref{D-3}) and (\ref{D-8}), we have
\begin{equation}\label{D-12}
\lim\limits_{k\rightarrow\infty}\int _{\Omega_k}e_0(\tilde{f}_{n,k})\,dA_0=m_i
\end{equation}
and
\begin{equation}\label{D-13}
\lim\limits_{k\rightarrow \infty}\int_{\Omega_k\setminus S^2_{+}}e_0(\tilde{f}_{n,k})\,dA_0=C_0,
\end{equation}
where $S^2_{+}$ is the northern hemisphere of $S^2$.

Choosing a sequence of compact sets $\{G_j\}$ such that $D_j\subset G_j \subset D_{j+1}$. Applying the Theorem \ref{Theorem 4.1} 
to the sequence $\tilde{f}_{n,k}$ on $G_j$ successively for $j$ and a diagonal argument, we obtain a subsequece of 
$\{\tilde{f}_{n,k}\}$ (still denoted by $\{\tilde{f}_{n,k}\}$), a finite subset $\{p_{i1},\ldots, p_{i\kappa_i}\}\subset S^2$ with $\kappa_i\leq \mathcal{A}_0/\epsilon_1$, 
and a Chern-harmonic map $f_{p_i}:S^2\longrightarrow M$ so that the subsequence converges to $f_{p_i}$ in $C^\infty$-topology on 
$S^2\setminus \{p_{i1},\ldots, p_{i\kappa_i},\textbf{\emph{s}}\}$, with the measures convergence
\begin{equation}\label{D-14}
e_0(\tilde{f}_{n,k})\, dA_0\rightarrow e_0(f_{p_i})\,dA_0+\sum\limits_{j=1}^{\kappa_i}m_{ij}\,\delta_{p_{ij}}+m_{i\textbf{\emph{s}}}\,\delta_{\textbf{\emph{s}}},
\end{equation}
where $m_{ij}>\epsilon_1$.  
This together with (\ref{D-13}) imply that there is no bubble points $p_{ij}$ in the southern hemisphere.
As the choosing of sequence $\{f_{n,k}\}$, passing to a subsequence, for any $j\geq k$, we can further assume that 
\begin{equation}\label{D-15}
\sup\limits_{x\in \partial D_k} \mbox{dist}\big(\tilde{f}_{n,j}(x), f_{p_i}(x)\big)\leq\frac{1}{k}
\end{equation}
and
\begin{equation}\label{D-16}
\sup\limits_{D'_k\setminus D_k}\big|e_0(\tilde{f}_{n,j})-e_0(f_{p_i})\big|\leq 1,
\end{equation}
where $D'_k:=\pi^{-1}\circ T_k\big(B_{4k\mu_k}(\tilde{x}_k)\big)\subseteq S^2$. Similarly, as the proof of (\ref{D-10}) and (\ref{D-11}), these imply
\begin{equation}\label{D-17}
\sup\limits_{S^2}\big|df_{p_i}\big|^2\leq C_3'\,E_0(f_{p_i}),
\end{equation}
and for any $j\geq k$ that
\begin{equation}\label{D-18}
\mbox{dist}\big(\tilde{f}_{n,j}(\partial D_k), f_{p_i}(\textbf{\emph{s}})\big)\leq \frac{C_4'}{k},
\end{equation}
where $C_3'$ depends on $(S^2, ds^2_0)$, $C_4'$ depends on $C_3'$ and $E_0(f_{p_i})$. Notice that $\partial B_{r_k}(\tilde{x}_k)\subseteq B_{2r_1}(0)\setminus B_{r_k/16k^2}(0)$,
so for any $x\in\partial B_{r_k}(\tilde{x}_k)$, $y\in\partial B_{2r_k}(0)$, (\ref{D-5}) and (\ref{D-11}) give 
\begin{eqnarray*}
\mbox{dist}\big(f_{n,k}(x),f_0(p_i)\big)&\leq&\mbox{dist}\big(f_{n,k}(x),f_{n,k}(y)\big)+\mbox{dist}\big(f_{n,k}(y),f_0(p_i)\big)\\
&\leq& 4r_k\;\sup\limits_{B_{2r_1}(0)\setminus B_{r_k/16k^2}(0)}\;\big|df_{n,k}\big|^2+\frac{C_2'}{k}\\
&\leq&\frac{C_5'}{k},
\end{eqnarray*}
where $C_5'$ depends on $C_2'$ and $\sup\limits_{\Sigma} e_0(f_0)$. This implies
\begin{equation}\label{D-19}
\mbox{dist}\big(f_{n,k}(\partial B_{r_k}(\tilde{x}_k)),f_0(p_i)\big)\leq \frac{C_5'}{k}. 
\end{equation}

We next clearly describe the position of $f_0(\Sigma)$ and $f_{p_i}(S^2)$ lie in $M$. We divide the energy on $B_4(0)$ into three parts as follows
\begin{eqnarray}\label{D-20}
\int_{B_4(0)}e_0(f_{n,k})\,dA_0&=&\int_{B_4(0)\setminus B_{r_k}(\tilde{x}_k)}e_0(f_{n,k})\,dA_0+\int_{B_{r_k}(\tilde{x}_k)\setminus B_{k\mu_k}(\tilde{x}_k)}e_0(f_{n,k})\,dA_0\nonumber\\
&{}&+\int_{B_{k\mu_k}(\tilde{x}_k)}e_0(f_{n,k})\,dA_0.
\end{eqnarray}
So, it is natural to define the \emph{neck map} 
\begin{equation}\label{D-21}
f_{n,k}|_{A_k}:A_k\longrightarrow M,\hspace{0.5cm}A_k:=B_{r_k}(\tilde{x}_k)\setminus B_{k\mu_k}(\tilde{x}_k),
\end{equation}
and the \emph{bubble map}
\begin{equation}\label{D-22}
\tilde{f}_{n,k}:D_k\longrightarrow M,
\end{equation}
whose image are agree with $f_{n,k}(B_{k\mu_k}(\tilde{x}_k))$. The corresponding domains $A_k$, $B_{k\mu_k}(\tilde{x}_k)$ will be called the \emph{neck domain}, \emph{bubble domain}, respectively. 
Notice that $f_{n,k}$ converges to $f_0$ in $C^\infty$-topology on $\Sigma\setminus \{p_1, \ldots, p_\kappa\}$, we have
\begin{equation}\label{D-23}
\lim\limits_{k\rightarrow\infty}\int_{B_4(0)\setminus B_{r_k}(\tilde{x}_k)}e_0(f_{n,k})\,dA_0=\int_{B_4(0)} e_0(f_0)\,dA_0.
\end{equation}
If we can show
\begin{equation}\label{D-24}
\lim\limits_{k\rightarrow\infty}\int_{B_{r_k}(\tilde{x}_k)\setminus B_{k\mu_k}(\tilde{x}_k)}e_0(f_{n,k})\,dA_0=0,
\end{equation}
then (\ref{C-5}), (\ref{D-20}) and (\ref{D-24}) implies
\begin{equation}\label{D-25}
\lim\limits_{k\rightarrow\infty}\int_{B_{k\mu_k}(\tilde{x}_k)}e_0(f_{n,k})\,dA_0=m_i.
\end{equation}

To prove (\ref{D-24}), we reparameterize the neck domain $A_k$ in polar coordinate as
\begin{equation*}
\phi_k:[0, T_k]\times S^1\longrightarrow \overline{A_k},\;(t,\theta)\mapsto r_ke^{-t+\texttt{i}\theta},
\end{equation*}
where $T_k=\log(r_k/k\mu_k)\rightarrow \infty$ as $k\rightarrow \infty$. We define a family of loops by
$\gamma_{k,t}:=f_{n,k}\circ\phi(t,\cdot)$, $t\in [0, T_k]$.

\begin{Proposition}\label{Proposition 5.1}
For the length of loops $\gamma_{k,t}$ with $t\in[0, T_k]$, for sufficient large $k$, we have
\begin{equation*}
\mathcal{L}(\gamma_{k,t})\leq 2\pi\sqrt{C_0C_3C'_6},
\end{equation*}
where $C_6'$ is a constant dpends on $ds^2_0$. Moreover, we have
\begin{equation*}
\mbox{\emph{max}}\big\{\mathcal{L}(\gamma_{k,0}), \mathcal{L}(\gamma_{k,T_k})\big\}\leq \frac{C_6}{k},
\end{equation*}
where $C_6$ is a constant depends on the metric $ds^2_0$ and $f_0$.
\end{Proposition}

\emph{Proof}. By the definition of $\mu_k$ in (\ref{D-8}), we have $E_0(f_{n,k}|_{A_k})\leq C_0$. So, by the Theorem \ref{Theorem 2.4}, on the unit disk around each point in $[1,T_k-1]\times S^1$,
we have $e_0(f_{n,k}|_{A_k})\leq C_0C_3$. For $1\leq t \leq T_k-1$, as the estimation given in (\ref{B-13}), we have
\begin{eqnarray*}
\mathcal{L}^2(\gamma_{k,t})&=&\big(\int_0^{2\pi}\big|df_{n,k}(\partial_\theta)\big|\;\big)^2\\
&\leq& C_6'\big(\int_0^{2\pi}d\theta\;\big)\;\big(\int_0^{2\pi}\big|df_{n,k}\big|^2d\theta\;\big)\\
&\leq& 4\pi^2C_0C_3C_6',
\end{eqnarray*}
where $C_6'$ is a constant depends on $ds^2_0$. For $0\leq t\leq 1$, since $\gamma_{k,t}\subseteq f_{n,k}(\overline{B_{r_k}(\tilde{x}_k)}\setminus B_{r_k/e}(\tilde{x}_k))$
and $\overline{B_{r_k}(\tilde{x}_k)}\setminus B_{r_k/e}(\tilde{x}_k)\subset B_{2r_1}(0)\setminus B_{r_k/16k^2}(0)$ for sufficient large k, 
it follows from (\ref{D-5}) and (\ref{D-10}) we have
\begin{eqnarray*}
\sup\limits_\theta\big|d\gamma_{k,t}(\partial_\theta)\big|^2&\leq&C_6'\sup\limits_{B_{r_k}(\tilde{x}_k)\setminus B_{r_k/e}(\tilde{x}_k)}r^2\,\big|df_{n,k}\big|^2\\
&\leq& \frac{C_6'C_7'}{k^2}, 
\end{eqnarray*}
where $C_7'$ depends on $f_0$. This implies that $\mathcal{L}(\gamma_{k,0})\leq C_6/k$, where $C_6$ depends on $C'_6$ and $C'_7$. For $T_k-1\leq t \leq T_k$, a similar estimation follows 
from (\ref{D-16}) and (\ref{D-17}).

\hfill$\Box$

\begin{Proposition}\label{Proposition 5.2}
For the constants $m_{i\textbf{s}}$ in $(5.14)$, we have
\begin{equation*}
m_{i\textbf{s}}=\limsup\limits_{k\rightarrow\infty} \int_{A_k} e_0(f_{n,k})\,dA_0.
\end{equation*}
\end{Proposition}

\emph{Proof}. Since $\Omega_k$ exhaust $S^2\setminus\{\emph{\textbf{s}}\}$ as $k$ goes to infinity, it follows from (\ref{D-14}) we obtain
\begin{eqnarray}\label{D-26}
m_{i\textbf{s}}&=&\lim\limits_{j\rightarrow\infty}\lim\limits_{k\rightarrow\infty} E_0(\tilde{f}_{n,k}|_{\Omega_k\setminus D_j})\nonumber\\
&=&\lim\limits_{j\rightarrow\infty}\lim\limits_{k\rightarrow\infty} \big( E_0(\tilde{f}_{n,k}|_{\Omega_k\setminus D_k})
+ E_0(\tilde{f}_{n,k}|_{D_k\setminus D_j})\big)\nonumber\\
&=&\lim\limits_{j\rightarrow\infty}\lim\limits_{k\rightarrow\infty}\big(E_0(f_{n,k}|_{B_{2r_k}(0)\setminus B_{r_k}(\tilde{x}_k)})
+E_0(f_{n,k}|_{A_k})+E_0(\tilde{f}_{n,k}|_{D_k\setminus D_j})\big).
\end{eqnarray}
Since $B_{2r_k}(0)\setminus B_{r_k}(\tilde{x}_k)\subseteq B_{2r_k}(0)\setminus B_{r_k/16k^2}(0)$, then (\ref{D-3}) implies that the first term in
(\ref{D-26}) is equal to zero. On the other hand, for any $k\geq j$, the (\ref{D-16}) and (\ref{D-17}) give
\begin{equation*}
\sup\limits_{D_{j+1}\setminus D_j}\,e_0(\tilde{f}_{n,k})\leq \sup\limits_{D_{j+1}\setminus D_j}\,\big(\big|e_0(\tilde{f}_{n,k})-e_0(f_{p_i})\big|+|e_0(f_{p_i})|\big)\leq C_8',
\end{equation*} 
where $C_8'$ depends on $f_{p_i}$. This implies that $\sup\limits_{D_k\setminus D_\ell} e_0(\tilde{f}_{n_k})\leq C_8'$, 
and hence the third term in (\ref{D-26}) is also equal to zero. So, the statement holds by
taking a subsequence  converges to $\limsup\limits_{k\rightarrow\infty} E_0(f_{n,k}|_{A_k})$.

\hfill$\Box$

To be continue, we define a new map
\begin{equation*}
F_k:B_4(0)\longrightarrow M\times \mathbb{C},\;\;x\mapsto (f_{n,k}, x).
\end{equation*}
Notice that $F_k$ is conformal, so it is easy to check that $F_k$ is also Chern-harmonic with respect to the product metric $h+dxd\bar{x}$ on $M\times \mathbb{C}$. 
We will use $F_k$ to prove the no energy loss and necklessness.
This approach has been used by R. Schoen \cite{[S-84]}, J. Jost \cite{[J-91]} and T.H. Parker \cite{[P-96]}.

\begin{Proposition}\label{Proposition 5.3}
For the neck maps $f_{n,k}|_{A_k}$, we have
\begin{equation*}
\limsup\limits_{k\rightarrow\infty} \int_{A_k} e_0(f_{n,k})\,dA_0=0\hspace{0.3cm} and \hspace{0.3cm}
\limsup\limits_{k\rightarrow\infty}\sup\limits_{x,x'\in A_k}\mbox{\emph{dist}}\big(f_{n,k}(x),f_{n,k}(x')\big)=0.
\end{equation*}
\end{Proposition}

\emph{Proof}. For the energy of $F_k$ on $A_k$, we have
\begin{eqnarray}\label{D-27}
\int_{A_k} e_0(F_k)\,dA_0=\int_{A_k}\big(e_0(f_{n,k})+1\big)\,dA_0
\leq C_0+C_9'r_k^2,
\end{eqnarray}
where $C_9'$ depends on $ds^2_0$. This implies $\mathcal{A}(F_k(A_k))\leq \epsilon_1$ for sufficient large $k$. 
For the length of $F_k(\partial B_r(\tilde{x}_k))$ with $k\mu_k\leq r\leq r_k$, we have
\begin{equation}\label{D-28}
\mathcal{L}(F_k(\partial B_r(\tilde{x}_k)))\leq \mathcal{L}(f_{n,k}(\partial B_r(\tilde{x}_k))) +2\pi r.
\end{equation}
Applying the isoperimetric inequality to $F_k(A_k)$  with sufficient large $k$, we have
\begin{equation}\label{D-29}
E_0(f_{n,k}|_{A_k})\leq E_0(F_k|_{A_k})= \mathcal{A}(F_k|_{A_k})\leq C_{10}'\, \mathcal{L}^2 (F_k|_{\partial A_k}),
\end{equation}
where $C_{10}'$ depends on $(M\times \mathbb{C}, h+dxd\bar{x})$.
So, it follows from (\ref{D-28}), (\ref{D-29}) and Proposition \ref{Proposition 5.1} we obtain $E_0(f_{n,k}|_{A_k})\rightarrow 0$ as $k\rightarrow \infty$. 

\vspace{0.1cm}

To prove the second identity, we only to show that for any fixed $\epsilon >0$, there exists a point $z\in M\times \mathbb{C}$ such that $F_k(A_k)\subseteq B_\epsilon (z)$ for large $k$.
It follows from (\ref{D-28}) that there exist $z,\;z'\in M\times \mathbb{C}$ such that $F_k(\partial A_k)\subseteq B_{\epsilon/4}(z)\cup B_{\epsilon/4}(z')$ for sufficient large $k$.
We claim that $F_k(A_k)\subseteq B_{\epsilon/2}(z)\cup B_{\epsilon/2}(z')$. If not, for some large $k$, there is a point $y\in F_k(A_k)$, but $y\notin B_{\epsilon/2}(z)\cup B_{\epsilon/2}(z')$.
Then, the Corollary \ref{Corollary 3.2} gives
\begin{equation*}
\frac{1}{4C_{10}'}\epsilon^2\leq\mathcal{A}(F_k(A_k)\cap B_\epsilon(y)),
\end{equation*}
which is a contradiction with (\ref{D-29}) for sufficient large $k$. Notice that $F_k(A_k)$ is connected, we obtain that $F_k(A_k)\subseteq B_\epsilon(z)$ or $B_\epsilon(z')$ for large $k$.

\hfill $\Box$

\emph{Remark}. It follows from Proposition \ref{Proposition 5.2} and Proposition \ref{Proposition 5.3} that $m_{i\textbf{\emph{s}}}=0$, and hence (\ref{D-14}) reduces to
\begin{equation}\label{D-30}
e_0(\tilde{f}_{n,k})\, dA_0\longrightarrow e_0(f_{p_i})\,dA_0+\sum\limits_{j=1}^{\kappa_i}m_{ij}\,\delta_{p_{ij}}.
\end{equation}

We finally iterate the previous procedure to get bubble tree convergence for Chern-minimal immersions.
Recall that for a given sequence of conformal Chern-minimal immersions $\{f_n\}$ from Riemann surfaces
$(\Sigma, \emph{\texttt{j}})$ into compact Hermitian surface $(M, J, h)$ with the areas are uniformly bounded by $\mathcal{A}_0$, 
then the Sacks-Uhlenbeck's procedure say that there are a subsequence $\{f_{n,k}\}$, a Chern-minimal map $f_0:\Sigma\longrightarrow M$, 
a finite set of bubble points $\{p_1,\ldots, p_\kappa\}\subset \Sigma$ with concentrated energy $m_1,\ldots, m_\kappa>\epsilon_1$ so that $f_{n,k}$ 
converges to $f_0$ in $C^\infty$-topology on $\Sigma\setminus\{p_1,\ldots,p_\kappa\}$ and the energy convergence
\begin{equation}\label{D-31}
\lim\limits_{k\rightarrow\infty} E(f_{n,k})=E(f_0)+\sum\limits_{i_1=1}^\kappa m_{i_1}.
\end{equation}
The renormalization procedure around any bubble point $p_{i_1}\in \{p_1,\ldots, p_\kappa\}$ say that there are a sequence 
of renormalized Chern-minimal immersions $\tilde{f}_{n,k,i_1}$ from domains $\Omega_k\subset S^2$ exhaust $S^2\setminus\{\emph{\textbf{s}}\}$,
a Chern-minimal map $f_{p_{i_1}}:S^2\longrightarrow M$, a finite set of bubble points $\{p_{i_11},\ldots, p_{i_1\kappa_{i_1}}\}\subset S^2$  
with concentrated energy $m_{i_11},\ldots, m_{i_1\kappa_{i_1}}>\epsilon_1$ so that
$\tilde{f}_{n,k,i_1}$ 
converges to $f_{p_{i_1}}$ in $C^\infty$-topology on $S^2\setminus\{p_{i_11},\ldots, p_{i_1\kappa_{i_1}},\emph{\textbf{s}}\}$ and (\ref{D-12}) and (\ref{D-30}) imply the energy convergence
\begin{equation}\label{D-32}
m_{i_1}=\lim\limits_{k\rightarrow\infty} E(\tilde{f}_{n,k,i_1})=E(f_{p_{i_1}})+\sum\limits_{i_2=1}^{\kappa_{i_1}} m_{i_1i_2}.
\end{equation}
Sequentially, if we have obtained a Chern-minimal map $f_{p_{I'}}:S^2\longrightarrow M$ and
a finite set of bubble points $\{p_{I'1},\ldots,p_{I'k_{I'}}\}$ with concentrated energy $m_{I'1},\ldots,m_{I'\kappa_{I'}}>\epsilon_1$, by repeating the
renormalization procedure, around each bubble point $p_{I'i_\ell}\in \{p_{I'1},\ldots,p_{I'k_{I'}}\}$, there are a sequence 
of renormalized Chern-minimal immersions $\tilde{f}_{n,k,I'i_\ell}$ from domains $\Omega_k\subset S^2$ exhaust $S^2\setminus\{\emph{\textbf{s}}\}$,
a Chern-minimal map $f_{p_{I'i_\ell}}:S^2\longrightarrow M$, a finite set of bubble points $\{p_{I'i_\ell \;1},\ldots, p_{I'i_\ell\;\kappa_{I'i_\ell}}\}\subset S^2$  
with concentrated energy $m_{I'i_\ell\; 1},\ldots, m_{I'i_\ell\;\kappa_{I'i_\ell}}>\epsilon_1$ so that
$\tilde{f}_{n,k,I'i_\ell}$ 
converges to $f_{p_{I'i_\ell}}$ in $C^\infty$-topology on $S^2\setminus\{p_{I'i_\ell\; 1},\ldots, p_{I'i_\ell\;\kappa_{I'i_\ell}},\emph{\textbf{s}}\}$ and the energy convergence
\begin{equation}\label{D-33}
m_{I'i_\ell}=\lim\limits_{k\rightarrow\infty} E(\tilde{f}_{n,k,I'i_\ell})=E(f_{p_{I'i_\ell}})+\sum\limits_{i_{\ell+1}=1}^{\kappa_{I'i_\ell}} m_{I'i_\ell i_{\ell+1}},
\end{equation}
where $I'=i_1\cdots i_{\ell-1}$ is a multiplicity index in $\textbf{N}^{\ell-1}$. We will agree on the convention that $f_{p_0}$ stands for the Chern-minimal map $f_0$ from
$\Sigma$ into $M$.
Notice that (\ref{D-13}) and (\ref{D-30}) imply that $E(f_{p_I})\geq C_0$,
thus this procedure will terminate in finite steps less than $\mathcal{A}_0/C_0$.

\begin{Theorem}\label{Theorem 5.4}
Let $\{f_n\}$ be a sequence of conformal Chern-minimal immersions from Riemann surfaces
$(\Sigma, \texttt{j})$ into compact Hermitian surface $(M, J, h)$ with the areas are uniformly bounded by $\mathcal{A}_0$.
Then there are a subsequence $\{f_{n,k}\}$, a Chern-minimal immersion $f_0:\Sigma \longrightarrow M$, a finite set of renormalized 
Chern-minimal sequences $\{\tilde{f}_{n,k,I}\}$ and
a finite set of Chern-minimal two-spheres $f_{p_I}:S^2\longrightarrow M$ so that

$(1)$ The sequences $\{f_{n,k}\}$, $\{\tilde{f}_{n,k,I}\}$ converge to $f_0$, $f_{p_I}$ in $C^\infty$-topology on $\Sigma\setminus\{p_1,\ldots,p_\kappa\}$, 
$S^2\setminus\{p_{I1},\ldots, p_{I\kappa_I}, \textbf{s}\}$, respectively.

\vspace{0.1cm}

$(2)$ There is no energy loss. That is
\begin{equation}\label{D-34}
\lim\limits_{k\rightarrow\infty}E(f_{n,k})=E(f_0)+\sum\limits_{I}E(f_{p_I}).
\end{equation}

$(3)$ There is no distance bubbling. Namely, for each bubble point $p_I$, we have $f_{p_I}(\textbf{s})=f_{p_{I'}}(p_I)$ 
with indices $I'=i_1\cdots i_{\ell-1}$ and $I=i_1\cdots i_{\ell-1}i_\ell$.
\end{Theorem}

\emph{Proof}. The statement of no energy loss follows from (\ref{D-31})-(\ref{D-33}), and the no distance bubbling
follows from (\ref{D-18}), (\ref{D-19}) and Proposition \ref{Proposition 5.3}.

\hfill$\Box$. 

The consequence of Theorem \ref{Theorem 5.4} imply that the images of $\{f_{n,k}\}$ pointwisely converges to the connected image of $\{f_0, f_{p_I}\}$.
The same as T.H. Parker's proof (see Corollary 2.3 in \cite{[P-96]}) of harmonic maps, the no distance bubbling implies that the bubble tree limit preserves the homotopy class.
\begin{Corollary}\label{Corollary 5.4}
Let $\{f_n\}$ be a sequence of conformal Chern-minimal immersions from Riemann surfaces
$(\Sigma, \texttt{j})$ into compact Hermitian surface $(M, J, h)$ with the areas are uniformly bounded by $\mathcal{A}_0$. If each $f_n$ represents the same
homotopy class $\alpha$ in the set $[\Sigma, M]$ of free homotopy class, then
\begin{equation*}
\alpha=[f_0]+\sum\limits_{I}[f_{p_I}].
\end{equation*}
\end{Corollary}

\vspace{0.6cm}

\noindent \textbf{Acknowledgments}. 
This project is supported by the NSFC (No.11871445), the Stable Support for Youth Team in Basic Research Field, CAS(YSBR-001) and the Fundamental 
Research Funds for the Central Universities.


\begin{thebibliography}{2}

\bibitem{[CT-97]} J.Y. Chen and G. Tian, {\it Minimal surfaces in Riemannian $4$-manifolds}, Geom. Funct. Anal., 7(1997), 873-916.

\bibitem{[GT-83]} D. Gilbarg and N. Trudinger, {\it Elliptic Partial Differential Equations of Second Order} (2nd ed.), New York: Springer-Verlag (1983). 

\bibitem{[HL-10]} X.L. Han and J.Y. Li, {\it Symplectic critical surfaces in Kähler surfaces}, J. Eur. Math. Soc., 12(2010), 505-527. 

\bibitem{[HS-74]} D. Hoffman and J. Spruck, {\it Sobolev and isoperimetric inequalities for Riemannian submanifolds}, Comm. Pure Appl. Math., 27(1974), 715-727.

\bibitem{[J-91]} J. Jost, {\it Two-dimensional geometric variational problems}, Wiley, New York, 1991.

\bibitem{[J-11]} J. Jost, {\it Riemannian geometry and geometric analysis}, (6th ed.), Springer-Verlag (2011).

\bibitem{[LW-08]} F.H. Ling and Ch.Y. Wang, {\it The analysis of harmonic maps and their heat flows}, World Scientific (2008).

\bibitem{[M-66]} C.B. Morrey, {\it Multiple Integrals in the  Calculus of Variations}, New York: Springer-Verlag (1966).

\bibitem{[PW-93]} T.H. Parker and J.G. Wolfson, {\it Pseudo-holomorphic maps and bubble trees}, J. Geom. Anal., 3(1993), 63-97. 

\bibitem{[P-96]} T.H. Parker, {\it Bubble tree convergence for harmonic maps}, J. Diff. Geom., 44(1996), 595-633.

\bibitem{[PX-21]} C.K. Peng and X.W. Xu, {\it On the Chern-minimal surfaces in Hermitian surface}, arXiv:2112.02304.

\bibitem{[S-84]} R. Schoen, {\it Analytic aspects of the harmonic map problem}, Lectures on Partial Differential Equations (1984), S.S. Chern, ed., Berling: Springer, 321-358.

\bibitem{[SU-81]} J. Sacks and K. Uhlenbeck, {\it The existence of minimal immersions of $2$-spheres}, Ann. of Math., 113(1981), 1-24.

\bibitem{[W-84]} S. Webster, {\it Minimal surfaces in a K$\ddot{a}$hler surface}, J. Diff. Geom., 20(1984), 463-470.

\bibitem{[W-86]} S. Webster, {\it On the relation between Chern and Pontrjagin numbers}, Contemporary Math., No.49, Amer. Math. Soc., Providence, RI, 1986, 135-143.

\bibitem{[W-88]} J.G. Wolfson, {\it Gromov's compactness of pseudo-holomorphic maps and symplectic geometry}, J. Diff. Geom., 28(1989), 383-405.  

\bibitem{[W-89]} J.G. Wolfson, {\it Minimal surfaces in K$\ddot{a}$hler surfaces and Ricci curvature}, J. Diff. Geom., 29(1989), 281-294.

\bibitem{[Y-94]} R.G. Ye, {\it Gromov's compactness theorem for pseudo holomorphic curves}, Trans. Amer. Math. Soc., 342(1994), 671-694.

\end{thebibliography}

\end{document}